\documentclass[11pt]{amsart}
\usepackage{amssymb, latexsym}
\theoremstyle{plain}
\newtheorem{theorem}{Theorem}
\newtheorem{corollary}{Corollary}

\newtheorem*{2'}{Theorem 2'}
\newtheorem*{3'}{Theorem 3'}

\theoremstyle{remark}

\newtheorem*{Remark 1}{Remark 1}
\newtheorem*{Remark 2}{Remark 2}
\newtheorem*{Remark 3}{Remark 3}
\newtheorem*{Remark 4}{Remark 4}

\numberwithin{equation}{section}

\begin{document}

\title[ Free Boundary for Reaction-Diffusion  with Convection]
 {The Behavior of the  Free Boundary for  Reaction-Diffusion Equations with Convection in an Exterior Domain with Neumann or  Dirichlet Boundary Condition}

\author{Ross G. Pinsky}
\address{Department of Mathematics\\
Technion---Israel Institute of Technology\\
Haifa, 32000\\ Israel}
\email{ pinsky@math.technion.ac.il}
\urladdr{http://www.math.technion.ac.il/~pinsky/}

\subjclass[2000]{35R35, 35J61} \keywords{reaction-diffusion equation, free boundary, exterior domain,
Neumann boundary condition, boundary flux, Dirichlet boundary condition}
\date{}

\begin{abstract}
Let
\begin{equation*}
\mathcal{L}=A(r)\frac{d^2}{dr^2}-B(r)\frac d{dr}
\end{equation*}
be a second order elliptic operator and consider the reaction-diffusion equation  with Neumann boundary  condition,
\begin{equation*}
\begin{aligned}
&\mathcal{L}u=\Lambda u^p\ \text{for}\ r\in(R,\infty);\\
&u'(R)=-h;\\
&u\ge0\ \text{is minimal},
\end{aligned}
\end{equation*}
where $p\in(0,1)$, $R>0$, $h>0$ and $\Lambda=\Lambda(r)>0$.
This equation is the radially symmetric case of  an equation of the form
\begin{equation*}
\begin{aligned}
&Lu=\Lambda u^p\ \text{in}\ \mathbb{R}^d-D;\\
&\nabla u\cdot \bar n=-h\ \text{on}\ \partial D;\\
&u\ge0 \ \text{is minimal},
\end{aligned}
\end{equation*}
where
\begin{equation*}
L=\sum_{i,j=1}^da_{i,j}\frac{\partial^2}{\partial x_i\partial x_j}-\sum_{i=1}^db_i\frac{\partial}{\partial x_i}
\end{equation*}
is  a second order elliptic operator, and where
 $d\ge2$,  $h>0$ is continuous, $D\subset R^d$ is bounded, and  $\bar n$ is the unit inward normal to the domain $\mathbb{R}^d-\bar D$.
Consider also the same equations with the Neumann  boundary  condition replaced by the Dirichlet boundary condition; namely,
$u(R)=h$ in the radial case and $u=h$ on $\partial D$ in the general case.
The solutions to the above  equations  may possess a free boundary. In the  radially symmetric case, if
$r^*(h)=\inf\{r>R:u(r)=0\}<\infty$, we call this   the radius of the free boundary; otherwise there is no free boundary.  We normalize the diffusion coefficient $A$ to be on unit order,  consider
the convection vector field $B$ to be on order $r^m$, $m\in R$, pointing either inward $(-)$ or outward $(+)$, and consider the reaction coefficient $\Lambda$ to be on order $r^{-j}$, $j\in R$. For both the Neumann boundary  case and the Dirichlet boundary case, we show for which choices of $m$, $(\pm)$ and $j$ a free boundary exists, and when it exists,
we obtain its growth rate  in $h$ as a function of $m$, $(\pm)$ and $j$. These results are then used to study the free boundary in  the non-radially symmetric case.

\end{abstract}

\maketitle

\section{Introduction and Statement of  Results}

Let $D\subset \mathbb{R}^d$, $d\ge2$,  be a bounded open set with smooth boundary
 such that $\mathbb{R}^d-\bar D$ is connected.
Let
\begin{equation}\label{nondivL}
L=\sum_{i,j=1}^da_{i,j}\frac{\partial^2}{\partial x_i\partial x_j}-\sum_{i=1}^db_i\frac{\partial}{\partial x_i}
\end{equation}
be a  strictly elliptic operator in $\mathbb{R}^d-D$ with smooth coefficients $a=\{a_{i,j}\}_{i,j=1}^d$ and $b=\{b_i\}_{i=1}^d$,
and let $\Lambda>0$ be a smooth function on $\mathbb{R}^d-D$.

In this paper, we consider two reaction-diffusion equations with sub-linear absorption, which are  of the same form,  but  which
have different boundary conditions, one the
Neumann condition   and the other the Dirichlet condition. We first describe the case of the Neumann boundary  condition.
We consider the
following reaction-diffusion equation in the exterior domain $\mathbb{R}^d-\bar D$,
\begin{equation}\label{nonLHE}
\begin{aligned}
&Lu=\Lambda u^p\ \text{in}\ \mathbb{R}^d-\bar D;\\
&\nabla u\cdot \bar n=-h\ \text{on}\ \partial D;\\
&u\ge0 \ \text{is minimal},
\end{aligned}
\end{equation}
where $p\in(0,1)$,  $h>0$ is continuous, and where $\bar n$ is the unit inward  to the domain $\mathbb{R}^d-\bar D$.
By minimal, we mean that  the solution  $u$ satisfies $u=\lim_{n\to\infty} u_n$, where
for large $n$, $u_n$ solves
\begin{equation}\label{nonLHE-n}
\begin{aligned}
&L u=\Lambda u^p \ \text{in}\ B_n-\bar D;\\
&\nabla u\cdot\bar n=-h\ \text{on}\ \partial D;\\
&u=0\ \text{on}\ \partial B_n,
\end{aligned}
\end{equation}
with $B_n$ denoting the ball of radius $n$ centered at the origin.
Existence for \eqref{nonLHE-n} will be shown in section \ref{existence}  by the method of upper and lower solutions. Uniqueness for \eqref{nonLHE-n},
the nonnegativity of $u_n$ and  the fact that $u_n$ is increasing in $n$ all follow
from a standard maximum principle argument for semi-linear equations.
From these facts, we obtain uniqueness for the solution to \eqref{nonLHE}.
The maximum principle shows that $u_n$ attains its maximum on $\partial D$.
Thus, we also obtain existence for \eqref{nonLHE} if the sequence
$\{u_n(x)\}_{n=1}^\infty$ is point-wise bounded for $x\in \partial D$.
The solution to \eqref{nonLHE} also attains its maximum on $\partial D$.

When $d=3$, $u$ can be thought of as the equilibrium quantity of a reactant  after having  undergone a long period of $L$-diffusion
and convection with sub-linear $p$-th power absorption with absorption coefficient $\Lambda$ in an exterior domain which is being supplied with the reactant via a normal boundary  flux $h$,  and where  complete and instantaneous   absorption
occurs far away. Note that the convection term in $L$ has been written as $-\sum_{i=1}^db_i\frac{\partial}{\partial x_i}$; under this
convention, the reactant is convected
in the direction $b$. We call $b$ the convection vector field.

A priori, it is not clear that a solution exists to \eqref{nonLHE} for every operator $L$ and every reaction coefficient $\Lambda$.
If the convection vector field $b$ points very strongly outward and/or if the reaction coefficient $\Lambda$ is very small, then one could imagine that the absorption term cannot overcome
the boundary flux and convection, leading to $\lim_{n\to\infty}u_n=\infty$.
(In the linear case, $\Lambda\equiv0$, a solution exists if and only if $L$ possesses an appropriate Green's function \cite{P14}.)
However, we believe that in fact the  solution always exists, as Theorem \ref{alwaysexist} below will  suggest.

We now turn to the corresponding reaction-diffusion equation with the Dirichlet boundary condition; namely, the equation
\begin{equation}\label{nonLHE-Dir}
\begin{aligned}
&Lv=\Lambda v^p\ \text{in}\ \mathbb{R}^d-\bar D;\\
&v=h\ \text{on}\ \partial D;\\
&v\ge0 \ \text{is minimal},
\end{aligned}
\end{equation}
where  $h>0$ is continuous. Similar to the previous case, minimality means that the solution $v$ satifies
$v=\lim_{n\to\infty}v_n$ where for large $n$, $v_n$ solves
\begin{equation}\label{nonLHE-n-Dir}
\begin{aligned}
&Lv=\Lambda v^p\ \text{in}\ B_n^d-\bar D;\\
&v=h\ \text{on}\ \partial D;\\
&v=0\ \text{on}\ \partial B_n.\\
\end{aligned}
\end{equation}
Existence for \eqref{nonLHE-n-Dir} will be shown in section \ref{existence} by the method of upper and lower solutions.
 Uniqueness for \eqref{nonLHE-n-Dir},
the nonnegativity of $v_n$ and  the fact that $v_n$ is increasing in $n$ all follow
from a standard maximum principle argument for semi-linear equations.
From these facts, we obtain uniqueness for the solution to \eqref{nonLHE-Dir}.
The maximum principle shows that $v_n$ attains its maximum on $\partial D$. Since $v_n=h$ on $\partial D$,
we conclude that $\{v_n(x)\}_{n=1}^\infty$ is point-wise bounded for $x\in D$, and thus
$v(x)\equiv\lim_{n\to\infty}v_n(x)$ is the unique solution to
 \eqref{nonLHE-Dir}.

The interesting  phenomenon that arises in the case of  sub-linear absorption
 is the possibility of a free boundary.
Specifically, there may exist an open set $\Omega$ satisfying  $D\subset \Omega\subset \mathbb{R}^d$ and  such that the smooth solution to
 \eqref{nonLHE}  or to \eqref{nonLHE-Dir} is positive in $\Omega-D$ and
identically zero in $\mathbb{R}^d-\Omega$.
The reason that this can occur is that when the solution  gets small, the absorption
term  is still relatively large since $p\in(0,1)$.

The free boundary for
 equation \eqref{nonLHE-Dir} in the case that the domain is not an exterior domain, but rather
  a bounded domain  or   an infinite slab, has
 been investigated by numerous authors. In many of these papers, the boundary value $h$ is fixed at 1 (the solution is interpreted as
 a concentration level of a reactant),
 the absorption coefficient $\Lambda$ is a constant $\lambda$, and one
 studies the critical value $\lambda^*$ such that the solution $v=v_\lambda$ satisfies  $\inf v_\lambda>0$ if
$\lambda<\lambda^*$ and $\inf v_\lambda=0$ if $\lambda\ge\lambda^*$.
For $\lambda\ge\lambda^*$, the
 region where $v\equiv0$ is known
as the ``dead core.'' See, for example, \cite{BS85} \cite{P88}, \cite{P92}, \cite{S96}.

For equations \eqref{nonLHE} and \eqref{nonLHE-Dir},
the existence and the location of a free boundary depend on the operator $L$, on the reaction term $\Lambda$ and on the boundary term $h$.
In this paper we study this dependence.

In order to investigate the free boundary quantitatively, we will  first and foremost  consider the radially symmetric case.
This is the case in which  $D=B_R$,   $h$ is a constant,
 $\Lambda(x)$ depends only on $|x|$  and $L$ is of the form
\begin{equation*}
L=A(|x|)\Delta-\hat B(|x|)\thinspace \frac x{|x|}\cdot \nabla.
\end{equation*}
When $L$ is of the above form, we will call it a radially symmetric operator.
The results obtained for  the  radially symmetric case will be used to obtain   results  for the general case.
 We denote by  $\mathcal{L}$ the radial part of  $L$; that is
\begin{equation}\label{radsymop}
\mathcal{L}=A(r)\frac{d^2}{dr^2}-B(r)\frac d{dr},
\end{equation}
where
$B(r)=\hat B(r)-\frac{(d-1)A(r)}r$.
 In the radially symmetric  case,
by uniqueness, the solutions to \eqref{nonLHE} and \eqref{nonLHE-Dir} are radial, and thus \eqref{nonLHE} for a function
$u(x)$ and \eqref{nonLHE-Dir}
for a function $v(x)$ reduce to the following
equations for  functions $u(r)$ and $v(r)$ respectively, where $r=|x|$:
\begin{equation}\label{nonLHE-rad}
\begin{aligned}
&\mathcal{L}u=\Lambda u^p\ \text{for}\ r\in(R,\infty);\\
&u'(R)=-h;\\
&u\ge0\ \text{is minimal};
\end{aligned}
\end{equation}
\begin{equation}\label{nonLHE-rad-Dir}
\begin{aligned}
&\mathcal{L}v=\Lambda v^p\ \text{for}\ r\in(R,\infty);\\
&v(R)=h;\\
&v\ge0\ \text{is minimal},
\end{aligned}
\end{equation}
where $h>0$ is a constant, $\Lambda=\Lambda(r)$ and
$\mathcal{L}$ is as in \eqref{radsymop}.
We call $B$  the radial convection vector field.
As before, by minimal we  mean that $u=\lim_{n\to\infty}u_n$ and $v=\lim_{n\to\infty}v_n$, where for large $n$,
$u_n$ and solves
\begin{equation}\label{nonLHE-rad-n}
\begin{aligned}
&\mathcal{L}u=\Lambda u^p\ \text{for}\ r\in(R,n);\\
&u'(R)=-h;\\
&u(n)=0,
\end{aligned}
\end{equation}
and $v_n$ solves
\begin{equation}\label{nonLHE-rad-n-Dir}
\begin{aligned}
&\mathcal{L}v=\Lambda v^p\ \text{for}\ r\in(R,n);\\
&v(R)=h;\\
&v(n)=0.
\end{aligned}
\end{equation}


For the radially symmetric case,
we define
$$
r_N^*=\inf\{r\ge R: u(r)=0\},\ \ r_D^*=\inf\{r\ge R: v(r)=0\}.
$$
If $r_N^*<\infty$ ($r_D^*<\infty$), we call $r_N^*$ ($r_D^*$) the radius of the free boundary for \eqref{nonLHE-rad} (for  \eqref{nonLHE-rad-Dir}).
In this case,
$u(r)=0$ ($v(r)=0$), for all $r\ge r_N^*$ ($r\ge r_D^*$). If $r_N^*=\infty$ ($r_D^*=\infty$), then  there is no free boundary.
Of course, $r_N^*$ and $r_D^*$  depend on the radius $R$ of the open set $D=B_R$, but
$R$ is fixed throughout this paper so we suppress this dependence. We write $r_N^*=r_N^*(h)$ and $r_D^*=r_D^*(h)$
to denote the dependence of the radius of the free boundary on the Neumann or Dirichlet boundary
 value $h$.
Under the assumption
that $A,B$ and $\Lambda$ have  power order  growth or decay as $r\to\infty$,
we will investigate whether or not $r_N^*(h)$ and $r^*_D(h)$  are finite, and in the case that  they are finite,
we will investigate how the asymptotic behavior of $r_N^*(h)$ and of $r_D^*(h)$ as $h\to\infty$ depend on these power orders.
Dividing \eqref{nonLHE-rad} by
the power order of $A$ allows us to normalize. Thus, in the sequel we will assume that $A$ is bounded and bounded from 0:
\begin{equation}\label{unifellrad}
C_1\le A(r)\le C_2,\     \text{for}\ r\ge R\ \text{where}\ 0<C_1<C_2.
\end{equation}
It turns out that the asymptotic behavior of $r_N^*(h)$, and to a bit lesser degree, of
$r_D^*(h)$, depend in a quite interesting and complicated fashion on
the sign and power order of $B$ and on the power order of $\Lambda$.

Before investigating the behavior of $r_N^*(h)$ and $r_D^*(h)$,
we present  a result which suggests that a solution
to \eqref{nonLHE-rad} exists for every operator $\mathcal{L}$ as in \eqref{radsymop}, regardless of how strongly outward the convection vector
field $B$ might point or how small the reaction coefficient $\Lambda$ might be.

Let $\exp^{(n)}(x)$ denote the $n$th iterate of  $e^x$; that is, $\exp^{(1)}(x)=e^x$ and $\exp^{(n)}(x)=e^{\exp^{(n-1)}(x)}$, $n\ge2$.
\begin{theorem}\label{alwaysexist}
 Assume that
$\mathcal{L}$ is as in \eqref{radsymop} with $A$ satisfying \eqref{unifellrad}.
Assume that for some positive integer $N$, one has
$B(x)\le\exp^{(N)}(x)$ and
 $\Lambda(r)\ge(\exp^{(N)}(r))^{-1}$.
 Then the solution to \eqref{nonLHE-rad} exists.
\end{theorem}

We now turn to our main focus,
the question of whether $r_N^*(h)$ and $r_D^*(h)$  are finite, and if they are,
the asymptotic behavior of $r_N^*(h)$ and $r_N^*(h)$  as $h\to\infty$.
We begin with the case in which the solutions to \eqref{nonLHE-rad}
and \eqref{nonLHE-rad-Dir} can be specified explicitly; namely the case that
$A$ and $\Lambda$ are  constant and $B\equiv0$.
One looks for the solution $u$ to \eqref{nonLHE-rad}
in the form $u(r)=\gamma(c-r)^l$. In order that this solve the differential
equation one is led to $l=\frac2{1-p}$ and $\gamma=\gamma_{\Lambda,A,p}\equiv\big(\frac{\Lambda(1-p)^2}{2A(1+p)}\big)^{\frac1{1-p}}$.
Then in order that $u(r)=\gamma_{\Lambda,A,p}(c-r)^{\frac2{1-p}}$ satisfy  the Neumann boundary condition, one is led to
$c=R+\big(\frac{(1-p)h}{2\gamma_{\Lambda,A,p}}\big)^{\frac{1-p}{1+p}}$. Thus, the explicit solution is
$$
u(r)=\begin{cases}\big(\frac{\Lambda(1-p)^2}{2A(1-p)}\big)^{\frac1{1-p}}\Big(r_N^*(h)-r\Big)^\frac2{1-p},& \ R\le r< r_N^*(h);\\  0, & \ r\ge r_N^*(h),
\end{cases}
$$
where the radius of the free boundary $r_N^*(h)$ is given by
\begin{equation}\label{explicit}
r_N^*(h)=R+\frac1{1-p}\big(\frac{2^p(1+p)A}\Lambda\big)^\frac1{1+p}\thinspace h^\frac{1-p}{1+p}.
\end{equation}
A similar calculation for \eqref{nonLHE-rad-Dir} gives
$$
v(r)=\begin{cases}\big(\frac{\Lambda(1-p)^2}{2A(1-p)}\big)^{\frac1{1-p}}\Big(r_D^*(h)-r\Big)^\frac  2{1-p},& \ R\le r< r_D^*(h);\\  0, & \ r\ge r_D^*(h),
\end{cases}
$$
where the radius of the free boundary $r_D^*(h)$ is given by
$$
r_D^*(h)=R+\frac1{1-p}\big(\frac{2(1+p)A}\Lambda\big)^\frac12h^\frac{1-p}2.
$$

Note that when $p\to1$, we have $r_N^*(h),r_D^*(h)\to\infty$, as expected, since there cannot be a free boundary in the linear case.
Note also that the smaller $p$ is, the  larger the order
of the free boundary as a function of $h$. At first glance this might seem surprising since, for fixed $h$,  the radius of the free boundary approaches $\infty$ when
$p\to1$.  However there are two opposing phenomena at play in $p$-th order absorption.
Where the solution is very small, the $p$-th order absorption is stronger the smaller $p$ is, however
where the solution is large, the $p$-th order absorption is stronger the larger $p$ is. Since we are now considering large  $h$, which
causes the solutions in both the Neumann and the Dirichlet cases to take on  large values,
it is this latter phenomenon which causes the order  of the radius of the free boundary as $h\to\infty$ to be larger when $p$ is smaller.

When $h\to\infty$,  $r_N^*(h)$ grows on the order $h^\frac{1-p}{1+p}$ while
$r_D^*(h)$ grows on the order $h^\frac{1-p}2$.
Thus, for large $h$,  the free boundary forms  farther away in the case of the Neumann boundary condition than
in the case of the Dirichlet boundary condition. In the general case, this phenomenon persists when the convection vector field $B$
points outwards, but not necessarily when  $B$ points inwards.
See Remark 2 after Theorem \ref{-1} in particular, and also
 Remark 2 after Theorems \ref{outwards} and \ref{m>-1} and the remark after Theorem \ref{Laplacian}.

We now consider convection vector fields $B$ which point outward.
We begin with the case in which there is no free boundary for any $h$  and the borderline case where the existence of a free boundary depends
on the value of $h$.

\begin{theorem}\label{nofreebdry}
Consider the solutions $u$ to \eqref{nonLHE-rad} and $v$ to \eqref{nonLHE-rad-Dir}, where
$\mathcal{L}$ is as in \eqref{radsymop} with
 $A$ satisfying \eqref{unifellrad}. Let $r^*_N(h)$ and $r_D^*(h)$ denote the free boundary radii for $u$ and $v$ respectively.

 \noindent i. Assume that
 $$
 B\ge0
 $$
and that
 $$
\Lambda(r)\le \frac K{r^{2+\epsilon}},
$$
for some $K>0$ and some $\epsilon\ge0$. If $\epsilon>0$, then
 there is no free boundary for any $h>0$; that is $r_N^*(h)=r_D^*(h)=\infty$, for all $h>0$.
If $\epsilon=0$, then there is no free boundary for sufficiently large $h$ and there is a free boundary for sufficiently small $h$; that is
$r_N^*(h)=r_D^*(h)=\infty$, for sufficiently large $h$, and $r_N^*(h), r_D^*(h)<\infty$, for sufficiently small $h$.

\noindent ii. Assume that for some $K>0$,
$$
 B(r)\ge \frac{r^m}K,
$$
for some $m>-1$,
and
$$
\Lambda(r)\le \frac K{r^j},
$$
for some $j\le 2$.
If
$$
m+j>1,
$$
 then
 there is no free boundary for any $h>0$; that is $r_N^*(h)=r_D^*(h)=\infty$, for all $h>0$.
If
$$
m+j=1,
$$
 then
 there is no free boundary for sufficiently large  $h$ and there is a free  boundary for sufficiently small $h$; that is $r_N^*(h)=r_D^*(h)=\infty$, for sufficiently large  $h$, and  $r_N^*(h), r_D^*(h)<\infty$, for sufficiently small $h$.
\end{theorem}

\bf\noindent Remark.\rm\ It is interesting and a little surprising that there are parameter values for which the existence of a free
boundary depends on the boundary value  $h$.
\medskip

Still considering convection vector fields which point outward,  we now turn to the case where there is a free boundary for all $h$.
In the sequel we will write $f(r)\approx g(r)$ to indicate that there exist constants $c_1,c_2>0$ such that
$c_1g(r)\le f(r)\le c_2g(r)$, for all  $r$.

\begin{theorem}\label{outwards}
Consider the solutions $u$ to \eqref{nonLHE-rad} and $v$ to \eqref{nonLHE-rad-Dir}, where
$\mathcal{L}$ is as in \eqref{radsymop} with
 $A$ satisfying \eqref{unifellrad}. Let $r^*_N(h)$ and $r_D^*(h)$ denote the free boundary radii for $u$ and $v$ respectively.
Assume that
\begin{equation}\label{>-1}
B(r)\approx r^M,\
\text{for some}\ M>-1
\end{equation}
or that
 \begin{equation}\label{=-1}
0\le  B(r)=O(r^{-1}).
\end{equation}
Let
$$
m=\begin{cases} & M,\ \text{if}\ \eqref{>-1}\ \text{holds};\\ & -1,\ \text{if}\ \eqref{=-1}\ \text{holds}.\end{cases}
$$
Also assume that
$$
\Lambda(r)\approx r^{-j},
$$
for some $j<2$.

\noindent i. If
$$
m+j\in[0,1),
$$
then
$$
r_N^*(h)\approx h^{\frac{1-p}{(1-m-j)p}} \ \ \text{and}\ \  r_D^*(h)\approx h^{\frac{1-p}{1-m-j}};
$$
\noindent ii. If $m+j<0$, then
$$
r_N^*(h)\approx h^{\frac{1-p}{p-m-j}}\ \ \text{and}\ \ r_D^*(h)\approx h^{\frac{1-p}{1-m-j}}.
$$
\end{theorem}

\noindent \bf Remark 1.\rm\ The theorem shows that the quantity $m+j$ determines the order of the exponent in the radius of the free boundary
for both \eqref{nonLHE-rad} and \eqref{nonLHE-rad-Dir}.
For the reason noted in the paragraph after \eqref{explicit}, the
smaller $p$ is, the larger the order of the free boundary as a function of $h$. Interestingly, for $r^*_N$, when $p\to0$, this order
goes to $\infty$ if $m+j\ge0$ (the regime of less effective absorption), while this order
stays bounded if $m+j<0$ (the regime of more effective absorption). On the other hand, for $r_D^*$, the order always stays
bounded when $p\to0$.
\medskip

\noindent \bf Remark 2.\rm\ Note that $r_N^*$ growths faster than $r_D^*$.

\medskip

We now consider convection vector fields $B$ which point inward.
\begin{theorem}\label{m>-1}
Consider the solutions $u$ to \eqref{nonLHE-rad} and $v$ to \eqref{nonLHE-rad-Dir}, where
$\mathcal{L}$ is as in \eqref{radsymop} with
 $A$ satisfying \eqref{unifellrad}. Let $r^*_N(h)$ and $r_D^*(h)$ denote the free boundary radii for $u$ and $v$ respectively.
 Assume that  $B$ satisfies
$$
B(r)\approx -r^m,
$$
with $m>-1$.
Assume that for some $N>0$ one has
$$
r^{-N}\le \Lambda(r)\le r^N,\ \text{for large}\ r.
$$
Then
\begin{equation}\label{rlformula}
r_N^*(h),\thinspace r_D^*(h)\approx(\log h)^\frac1{1+m}.
\end{equation}
\end{theorem}

\noindent\bf Remark 1.\rm\ In contrast to the case in which the convection vector field points outward, when the convection
vector field points inward on a power order larger than $r^{-1}$, the order of the radius of the free boundary is
insensitive to the power order of the  reaction coefficient.

\medskip

\noindent \bf Remark 2.\rm\ Note that at logarithmic orders, $r_N^*$ growths at the same rate as $r_D^*$.
In light of Theorem \ref{-1} below, we suspect that in fact $r_D^*$ grows faster than $r_N^*$. See Remark
2 after Theorem \ref{-1}.
\medskip

It turns out that the asymptotic behavior of $r_N^*(h)$ is very sensitive to small changes in the convection vector field $B$ when
the convection vector field points inward and is on the order $\frac1r$. This sensitivity does not hold
for $r_D^*(h)$.
For this case, we present a result only in the case that the reaction coefficient is on unit order.
\begin{theorem}\label{-1}
Consider the solutions $u$ to \eqref{nonLHE-rad} and $v$ to \eqref{nonLHE-rad-Dir}, where
$\mathcal{L}$ is as in \eqref{radsymop} with $A=A_0$ and $B(r)=-\frac {B_0}r$, where
$A_0,B_0$ are positive constants. Let $r^*_N(h)$ and $r_D^*(h)$ denote the free boundary radii for $u$ and $v$ respectively.
Define
$$
\mu=\frac{B_0}{A_0}>0.
$$
Assume
$$
\Lambda\approx 1.
$$
Then
\begin{equation}\label{1rformula}
r_N^*(h)\approx h^\frac{1-p}{1+p+\mu(1-p)} \ \ \text{and}\ \ r_D^*(h)\approx h^\frac{1-p}2.
\end{equation}
\end{theorem}

\noindent \bf Remark 1.\rm\
Theorems \ref{-1} and \ref{outwards} demonstrate the particular sensitivity of $r_N^*(h)$ in the case that the convection
vector field is pointing inward  on the order $\frac1r$. Indeed,
when the reaction coefficient  $\Lambda$ is on unit order, the  diffusion coefficient is a constant and $B(r)=-\frac{B_0}r$,
Theorem \ref{-1} shows that the power order of $r_N^*(h)$
runs from $\frac{1-p}{1+p}$ to 0 as $B_0$ runs from 0 to $\infty$.
In contrast, Theorem \ref{outwards} shows that when $\Lambda$ and the diffusion coefficient are on unit order and
$B$ points outward on the order $\frac1r$, then $r_N^*(h)\approx h^{\frac{1-p}{1+p}}$, while if $B(r)\approx r^m$, for $m>-1$, then
the power order of
$r_N^*(h)$  runs from $\frac{1-p}{1+p}$ to $\infty$ as $m$ runs from $-1$ to 1.

\medskip

\noindent \bf Remark 2.\rm\ Note that $r_N^*$ growths faster than $r_D^*$ when $\mu\in(0,1)$, but
 $r_D^*$ growths faster than $r_N^*$ when $\mu>1$. When the
 convection vector field $B$ points outward, Theorem \ref{outwards} shows that
 $r_N^*$ growths faster than $r_D^*$, while when the convection vector field $B$ points inwards at a power order larger than $r^{-1}$, Theorem \ref{m>-1}
shows that  at logarithmic order, $r_N^*$ and $r_D^*$ grow at the same rate.
Now Theorem \ref{-1} shows that for a small band of inward pointing convection vector fields  $B$
whose power order is $r^{-1}$, that is,
whose strength and direction
fall between those considered in Theorems \ref{outwards} and \ref{m>-1},
$r_N^*$ grows faster for the ones that point inward more weakly  and $r_D^*$ grows faster for the ones
that point inward more strongly.
This would suggest that in fact $r_D^*$ grows faster than $r_N^*$  for the inward convection vector fields of Theorem \ref{m>-1}.

\medskip

We now consider the general non-radially symmetric case. Now $h$ is  a function defined on $\partial D$ instead of being a constant.
For our final two theorems below, we will write this function in  the form $hh_0$, where $h_0$ is a
continuous, strictly positive function on $\partial D$, and $h>0$ is a parameter.
For equations \eqref{nonLHE} and \eqref{nonLHE-Dir}, we
 define the \it free  boundary inner radii\rm\ $r_N^{*,-}(h)$ and $r_D^{*,-}$,  and the \it free boundary outer radii\rm\ $r_N^{*,+}(h)$
 and $r_D^{*,+}(h)$,  by
$$
\begin{aligned}
&r_N^{*,-}(h)=\inf\{|x|:u(x)=0\};\ \  r_N^{*,+}(h)=\sup\{|x|:u(x)>0\};\\
&r_D^{*,-}(h)=\inf\{|x|:v(x)=0\};\ \  r_D^{*,+}(h)=\sup\{|x|:v(x)>0\}.
\end{aligned}
$$
Recall that a subset $D\subset \mathbb{R}^d$ satisfying $0\in D$ is called
 \it star-shaped\rm\ with respect to 0 if the line segment connecting $0$ to $x$ is entirely contained in $D$, for all $x\in D$.

We begin by considering the $d$-dimensional Laplacian in a general  domain.
Using  Theorem \ref{-1}, we will prove the following result.
\begin{theorem}\label{Laplacian}
Consider the solutions $u$ to  \eqref{nonLHE} and $v$ to \eqref{nonLHE-Dir},  with $L=\Delta$ and $\Lambda\approx1$.
 In the case of \eqref{nonLHE},  assume that $D\subset R^d$ is star-shaped with respect to 0.
Let the boundary term in \eqref{nonLHE} and \eqref{nonLHE-Dir} be given by $hh_0$, where
$h_0$ is a continuous, strictly positive
function on $\partial D$, and $h>0$ is a parameter.
Then
$$
\begin{aligned}
&r_N^{*,+}(h)\approx r_N^{*,-}(h)\approx h^\frac{1-p}{1+p+(d-1)(1-p)};\\
&r_D^{*,+}(h)\approx r_D^{*,-}(h)\approx h^\frac{1-p}2.
\end{aligned}
$$
\end{theorem}

\noindent \bf Remark.\rm\  As an addendum to Remark 2 after Theorem \ref{-1},
note that in one dimension, $r_N^{*,\pm}$ grow faster than $r_D^{*,\pm}$,
in two dimensions, $r_N^{*,\pm}$ and  $r_D^{*,\pm}$ grow on the same power order, and
in three dimensions or higher,
$r_D^{*,\pm}$ grow faster than $r_N^{*,\pm}$.

\medskip

We now consider non-radially symmetric operators.
For an operator $L$ as in \eqref{nondivL}, defined on all of $\mathbb{R}^d$, define
\begin{equation}\label{radialb}
B^+(r)=\sup_{|x|=r}\frac xr\cdot b(x); \ \
B^-(r)=\inf_{|x|=r}\frac xr\cdot b(x),\ r>0.
\end{equation}
Analogous to \eqref{unifellrad}, we normalize by requiring the diffusion matrix $a$ to satisfy
\begin{equation}\label{unifell}
C_1\le \sum_{i,j=1}^da_{i,j}(x)v_iv_j\le C_2,   \ \text{for all}\ v\in \mathbb{R}^d\ \text{satisfying}\ |v|=1\ \text{and all}\ x\in \mathbb{R}^d.
\end{equation}
Define
\begin{equation}\label{AB}
AB^-(r)=\inf_{|x|=r}\Big(\sum_{i=1}^da_{i,i}(x)-\sum_{i,j=1}^da_{i,j}(x)\frac{x_i}r\frac{x_j}r +\sum_{i=1}^dx_ib_i(x) \Big).
\end{equation}
For the next theorem, we will use the following terminology.
If $f\approx g$,   we say that \it $F$ satisfies the upper (lower) bound satisfied by $f$\rm\ if
$F(r)\le cg(r)$ ($F(r)\ge cg(r)$), for all  $r$ and some $c>0$.

The theorem below converts the results of Theorems \ref{alwaysexist}-\ref{m>-1} to results for the non-radially symmetric case.
\begin{theorem}\label{nonradial}
Consider \eqref{nonLHE} and \eqref{nonLHE-Dir} with  $L$ as in \eqref{nondivL} defined on all of $\mathbb{R}^d$ and satisfying \eqref{unifell}.
Let the boundary term in \eqref{nonLHE} and \eqref{nonLHE-Dir} be given by $hh_0$, where
$h_0$ is a continuous, strictly positive
function on $\partial D$, and $h>0$ is a parameter.
Let $B^\pm$ be as in \eqref{radialb} and let $AB^-$ be as in \eqref{AB}.
For \eqref{nonLHE}, assume that $D$ is star-shaped with respect to 0.

\noindent i. Assume that $\Lambda$ satisfies the bound specified in Theorem \ref{alwaysexist} and that $B^+$ satisfies the bound
on $B$ specified in Theorem \ref{alwaysexist}. Then the solution to \eqref{nonLHE} exists;

\noindent ii.
Assume either that  $\Lambda$ satisfies the bound in Theorem \ref{nofreebdry}-i and  that    $AB^-\ge0$, or that  $\Lambda$ satisfies the bound
in Theorem  \ref{nofreebdry}-ii and $B^-$ satisfies the bound for $B$ in Theorem  \ref{nofreebdry}-ii. Then there is no free boundary
for \eqref{nonLHE} or \eqref{nonLHE-Dir};

\noindent iii.  Assume that $\Lambda$ satisfies the lower (upper) bound in Theorem \ref{outwards} and assume that
$B^+$ satisfies the upper  bound  ($B^-$ satisfies the lower bound)  for $B$ in \eqref{>-1} of  Theorem \ref{outwards}.
Then
the free boundary outer radii $r_N^{*,+}(h)$ and $r_D^{*,+}(h)$ (free boundary inner radii $r_N^{*,-}(h)$ and $r_D^{*,-}(h)$ )   satisfy the upper (lower) bounds satisfied by $r_N^*(h)$ and $r_D^*(h)$ in Theorem \ref{outwards};

\noindent iv. Assume that $\Lambda$ satisfies the lower (upper) bound in Theorem \ref{outwards}. Assume that
$B^+$ satisfies the upper bound for $B$ in \eqref{=-1} of  Theorem \ref{outwards} (that $AB^-$ satisfies $AB^-\ge0$).
Then
the free boundary outer radii $r_N^{*,+}(h)$ and $r_D^{*,+}(h)$ (free boundary inner radii $r_N^{*,-}(h)$ and $r_D^{*,-}(h)$)   satisfy the upper (lower) bounds satisfied by  $r_N^*(h)$ and $r_D^*(h)$ in Theorem \ref{outwards}
with $m=-1$;

\noindent v. Assume that  $\Lambda$ satisfies the bound in Theorem \ref{m>-1} and that $B^+$ satisfies the upper bound
($B^-$ satisfies the lower bound) for $B$ in Theorem \ref{m>-1}. Then the  free boundary outer radii $r_N^{*,+}(h)$ and $r_D^{*,+}(h)$
(free boundary inner radii $r_N^{*,-}(h)$ and $r_D^{*,-}(h)$)
satisfy the upper (lower) bounds satisfied by $r_N^*(h)$ and $r_D^*(h)$ in Theorem \ref{m>-1}.
\end{theorem}


We expect that the results in Theorems \ref{Laplacian} and  \ref{nonradial} for \eqref{nonLHE}
hold without the requirement that the open set $D$ be star-shaped. It was necessary for   our method of proof, which
 makes heavy use of point-wise comparison methods. Such methods
are more flexible with the Dirichlet boundary condition than they are with the Neumann boundary condition.

The method of proof used in all of the theorems   involves the construction of appropriate upper and lower solutions.
In certain instances, these  upper and lower solutions have a relatively  simple form, but frequently they are quite
complicated.

The rest of the paper is organized as follows. In section \ref{existence} we recall the theory of upper and lower solutions,
and deal with a technical issue that  arises  with regard to the construction of such solutions.
In section \ref{always}, we give a quick proof of Theorem \ref{alwaysexist} by constructing a  simple upper solution.
In section \ref{sec-nofree} we first give the rather delicate and involved proof of Theorem \ref{outwards}, by
constructing appropriate upper and lower solutions. Then we use
that proof to give a quick proof of Theorem \ref{nofreebdry}.
In section  \ref{logh} we construct fairly simple upper and lower solutions to prove Theorem \ref{m>-1}.
In section \ref{fuchsian}, we construct upper and lower solutions of a  different sort than
what has been constructed so far to prove Theorem \ref{-1}.
In section
\ref{generalcase} we prove
Theorems \ref{Laplacian} and \ref{nonradial} by constructing appropriate upper and lower solutions with the help of the
radially symmetric upper and lower solutions already constructed.


\section{Existence and Comparison Via the Method of Upper and Lower Solutions}\label{existence}

A $C^2$-function  $V_n^+\ge0$ is an  \it upper solution \rm\  to \eqref{nonLHE-n} if it satisfies $LV_n^+\le 0$ in $B_n-D$ and  $\nabla V_n^+\cdot\bar n\le -h$ on $\partial D$. A $C^2$-function $V_n^-\ge0$ is a \it lower solution\rm\ to \eqref{nonLHE-n} if it satisfies
 $LV_n^-\ge0$ in $B_n-D$, $\nabla V_n^-\cdot\bar n\ge -h$ on $\partial D$ and $V_n^-=0$ on $\partial B_n$.
Similarly, a
 $C^2$-function  $V_n^+\ge0$ is an  \it upper solution \rm\  to \eqref{nonLHE-n-Dir} if it satisfies $LV_n^+\le 0$ in $B_n-D$ and
  $V_n^+\ge h$ on $\partial D$, while a  $C^2$-function $V_n^-\ge0$ is a \it lower solution\rm\ to \eqref{nonLHE-n-Dir} if it satisfies
 $LV_n^-\ge0$ in $B_n-D$, $V_n^-\le h$ on $\partial D$ and $V_n^-=0$ on $\partial B_n$.
A standard application of the maximum principle for semi-linear equations shows that if there is a solution $u_n$ to \eqref{nonLHE-n} and $V_n^+$ is an upper solution
($V_n^-$ is a lower solution), then $u_n\le V_n^+$ ($u_n\ge V_n^-$). The corresponding statement also holds for
solutions $v_n$ to \eqref{nonLHE-n-Dir}.
A fundamental result states that if there exist an upper solution $V_n^+$ and a lower solution $V_n^-$ to \eqref{nonLHE-n} ( to \eqref{nonLHE-n-Dir}) such that
$V_n^-\le V_n^+$, then in fact there exists a solution $u_n$ to \eqref{nonLHE-n} ($v_n$ to \eqref{nonLHE-n-Dir}) \cite{S73}. The solution
$u_n$ ($v_n$)
necessarily satisfies $V_n^-\le u_n\le V_n^+$  ($V_n^-\le v_n\le V_n^+$).

Of course, the radially symmetric equations \eqref{nonLHE-rad-n}  and \eqref{nonLHE-rad-n-Dir}   are  particular cases of
\eqref{nonLHE-n} and \eqref{nonLHE-n-Dir} respectively.
In the sections that follow, we will construct radially symmetric upper and lower solutions $V_n^+, V_n^-$
to \eqref{nonLHE-rad-n} and \eqref{nonLHE-rad-n-Dir},
satisfying $V_n^-\le V_n^+$. This will prove
existence of the solution $u_n$ to \eqref{nonLHE-rad-n} and the solution $v_n$ to \eqref{nonLHE-rad-n-Dir} and also give upper and lower bounds on $u_n$ and $v_n$.
In fact, these upper and lower bounds will be independent of $n$ for sufficiently large $n$; thus they will also provide upper and lower
bounds for the solutions $u$ and $v$ to \eqref{nonLHE-rad} and \eqref{nonLHE-rad-Dir}.
In the last section, these radially symmetric upper and lower solutions are used to construct
appropriate upper and lower solutions for \eqref{nonLHE-n} and \eqref{nonLHE-n-Dir}, under the additional
assumption  that $D$ is star-shaped in the case of \eqref{nonLHE-n}. This will give existence and upper and lower bounds for the general case.

We deal here with one technical issue that will arise in the construction of the radially symmetric upper and lower solutions.
We will construct radially symmetric upper and lower solutions in the following form.
Let $0<R<c$, let $l>1$ and let $f^\pm:[R,c]\to (0,\infty)$ be  smooth functions.
Consider  the functions
$$
V_c^\pm(r)=\begin{cases} (c-r)^lf^\pm(x),& R\le r\le c,\\ 0,& r>c.\end{cases}
$$
Since $l>1$, it follows that $V$ is a $C^1$-function on $[r,\infty)$; however, since we have not assumed that $l>2$, $V$ is not necessarily
$C\thinspace^2$ at $r=c$.
We wish to show that the maximum principle for semi-linear equations still holds.
(We illustrate this with regard to the solution $u_n$ to \eqref{nonLHE-rad-n}; the same thing holds with regard to the solution $v_n$ to
\eqref{nonLHE-rad-n-Dir}.)
Namely, if $n>c$ and $u_n$ is  a solution to \eqref{nonLHE-rad-n},
we wish to show that if $V_c^-$ satisfies
$\mathcal{L}V_c^-\ge\Lambda (V_c^-)^p$, on $[R,c)$ and $(V_c^-)'(R)\ge-h$, then
$u_n(r)\ge V_c^-(r)$ on $[R,n]$, and if
$V_c^+$ satisfies $\mathcal{L}V_c^+\le(V_c^+)^p$ on $[R,c)$ and $(V_c^-)'(R)\le-h$, then $u_n(r)\le V_c^+(r)$ on $[R,n]$.
Since $V_c^\pm$ does not depend on $n$, letting $n\to\infty$ in the above inequalities
will then give
\begin{equation}\label{lowerbound}
u\ge V_c^-\ \text{on}\ [R,\infty),
\end{equation}
and
\begin{equation}\label{upperbound}
u\le V_c^+\ \text{on}\ [R,\infty),
\end{equation}
where $u$ is the solution to \eqref{nonLHE-rad}.
Now \eqref{lowerbound} implies that $r^*(h)\ge c$
while \eqref{upperbound} implies that $r^*(h)\le c$.

Nothing extra needs to be done in the first case. Indeed, since $V_c^-$ is a $C\thinspace^2$-function on $[R,c)$ and vanishes continuously at $c$,
it follows from the maximum
principle for semi-linear equations that $u_n\ge V_c^-$ on $[R,c]$. Then since $u_n\ge0$ on $[c,n]$ and $V_c^-\equiv0$ on $[c,n]$, we conclude that
$u_n\ge V_c^-$ on $[R,n]$.
Now consider the second case. Since $V_c^+-u_n$ is $C\thinspace^2$ on $[R,n]$ except at $c$, it
 follows by the maximum principle for semi-linear equations that either
$V_c^+-u_n\ge 0$ on $[R,n]$ or $V_c^+-u_n$ on $[R,n]$ attains a negative minimum  at $c$. We  show that this second
situation cannot occur.  For convenience, denote $V_c^+-u_n$ by $W$.
Note that $\lim_{r\to c^-}W''(r)$ exists as an extended real value; it is finite if $l=2$ and is equal to $+\infty$ if $l\in(1,2)$.
Assume now that $W$ attains a negative minimum at $c$.
It follows that
$W'(c)=0$ and $\lim_{r\to c^-}W''(r)\ge0$.
But this is impossible by the maximum principle for semi-linear equations. More precisely, an application of the mean value
theorem shows that $W$ solves an equation
of the form $\mathcal{L}W-HW=0$ in $[R,c)$, where $H=H(r)>0$ in $[R,c]$. ($H(c)>0$ because $u_n(c)>0$ since $W$  is assumed to be negative at $c$.) Thus the above conditions on $W$ yield
the contradiction $\lim_{r\to c^-}(\mathcal{L}W-HW)>0$.

\section{Proof of Theorem \ref{alwaysexist}}\label{always}
Let
$F(r)=\exp^{(N)}(r)$. So $B\le F$ and $\Lambda\ge\frac1{F}$. Let
$$
U(r)=\gamma+\frac h{F(R)}e^{-\int_R^rF(s)ds}.
$$
Note that $U'(R)=-h$. We will show that for an appropriate choice of $\gamma>0$, $U$ will satisfy
$\mathcal{L}U\le \Lambda U^p$. Thus, $U$ will be  an upper solution for \eqref{nonLHE-rad-n}, for all $n$. Since 0 is a lower solution, this  proves that a solution to
\eqref{nonLHE-rad} exists.

Recalling \eqref{unifellrad}, we
have
\begin{equation}\label{mathcalupper}
\mathcal{L}U(r)=\frac h{F(R)}\big(AF^2-AF'+BF)(r)e^{-\int_R^rF(s)ds}\le \frac{(C_2+1)h}{F(R)}F^2(r)e^{-\int_R^rF(s)ds}.
\end{equation}
We have
\begin{equation}\label{rhslower}
\Lambda(r) U^p(r)\ge \frac{\gamma^p}{F(r)}.
\end{equation}
Since $\sup_{r\ge R}F^3(r)e^{-\int_R^rF(s)ds}<\infty$, it follows from \eqref{mathcalupper}
and \eqref{rhslower} that if $\gamma$ is chosen sufficiently large, then $U$ is an upper solution.
\hfill $\square$
\medskip

\section{Proofs of Theorems \ref{nofreebdry} and \ref{outwards}}\label{sec-nofree}
We first prove Theorem \ref{outwards}.  Then we prove     Theorem \ref{nofreebdry}, which  will follow
 from the proof of Theorem \ref{outwards} after  making minor revisions.
\medskip

\noindent \it Proof of Theorem \ref{outwards}.\rm\
We present the proof for $r_N^*$. Upon its completion,  we explain how the proof for $r_D^*$ follows from the
proof for $r_N^*$. For the solution $u_n$ to \eqref{nonLHE-rad-n},
we look for upper and lower solutions
 in the form
\begin{equation}\label{V}
V(r)=\begin{cases}(c-r)^lf(r),&\ R\le r\le c\\ 0,&\ r>c,\end{cases}
\end{equation}
where $c>R$,
$l>1$ and $f>0$ is smooth.
From the discussion in section \ref{existence}, we may ignore the fact that $V$ might not be $C^2$ at $r=c$.
We emphasize that $f,c$ and $l$ will not depend on $n$, and thus $V$ will be independent of $n$.
The function $f$ and the value $c$ will depend on the boundary flux $h$.
If such a $V$ with $c=c(h)$  is an upper solution for all $n>c(h)$, then it follows that
the solution $u=\lim_{n\to\infty}u_n$ to \eqref{nonLHE-rad} satisfies  $u\le V$, and consequently,
the radius of the free boundary $r^*(h)$ satisfies $r^*(h)\le c(h)$. Similarly, if such a $V$ with $c=c(h)$ is a lower
solution, then it follows that $u\ge V$, and consequently,  $r^*(h)$ satisfies $r^*(h)\ge c(h)$.

In order that $V$ be an upper or lower solution to \eqref{nonLHE-rad-n}, it must satisfy a differential inequality in $[R,n]$ corresponding to the differential
equation, and a differential inequality at $R$ corresponding to the boundary condition.
Substituting $V$ into the differential equation in \eqref{nonLHE-rad-n}, it follows that in order that $V$ be an upper solution for all
large $n$, $f$ must satisfy the differential inequality
\begin{equation}\label{diff-ineq}
Af''-\big(\frac{2lA}{c-r}+B(r)\big)f'+\big(\frac{l(l-1)A}{(c-r)^2}+\frac{lB(r)}{c-r}\big)f\le\Lambda(r)(c-r)^{-l(1-p)}f^p, \ R\le r\le c,
\end{equation}
and in order that $V$ be a lower solution for all large $n$, $f$ must satisfy the above differential inequality with the direction of the inequality reversed.
Substituting  $V$ into the  boundary condition in \eqref{nonLHE-rad-n}, it follows that in order that $V$ be an upper solution, $f$ must satisfy
the boundary differential inequality
\begin{equation}\label{bdryinequ}
l(c-R)^{l-1}f(R)-(c-r)^lf'(R)\ge h,
\end{equation}
and that in order that $V$ be a lower solution, $f$ must satisfy the above boundary differential inequality with the direction of the inequality
reversed.
We are interested in  this inequality for large $h$.
As noted, we will consider $c$ to be a function of $h$.
The function $f$ will depend on $h$ only through its dependence on $c$.
The parameter $l$ will be chosen independent of $h$ and $c$.

We consider $f$ in the form
\begin{equation}\label{f}
f(r)=\theta c^{-\delta}(\gamma+r^{-k}),\ \text{where}\ \delta\in R,\ \gamma,\theta>0 \ \text{and}\
k=c^{-L},\ \text{with}\ L>0.
\end{equation}
The values $\theta,\gamma, \delta$ and $L$ will be chosen independent of $c$.
The inequality \eqref{bdryinequ} becomes
\begin{equation}\label{bdryinequagain}
\theta\big[l(c-R)^{l-1} c^{-\delta}(\gamma+R^{-k})+(c-R)^l c^{-\delta}kR^{-k-1}\big]\ge h.
\end{equation}
Thus, recalling that $k=c^{-L}$,  it follows that
if $\max(l-\delta-1,l-\delta-L)>0$, then
  condition \eqref{bdryinequagain} will be fulfilled for large $h$
if $c=c(h)\ge  Mh^{\frac1{\max(l-\delta-1,l-\delta-L)}}$, for some appropriate $M>0$ (which depends on the constants $l,R,\theta,\gamma$).
Similarly, \eqref{bdryinequagain} with the direction of the inequality reversed will be fulfilled if
$c=c(h)\le Mh^{\frac1{\max(l-\delta-1,l-\delta-L)}}$, for some appropriate $M>0$.
Thus, if we can find values $l^+,\delta^+$ and $L^+$ (as well as appropriate values for $\theta$ and $\gamma$)
such that \eqref{diff-ineq} holds, and if we also find values $l^-,\delta^-$ and $L^-$ (as well as appropriate values
of $\theta$ and $\gamma$) such that \eqref{diff-ineq} holds with the direction of the inequality reversed, and if furthermore,
$\max(l^+-\delta^+-1,l^+-\delta^+-L^+)=\max(l^--\delta^--1,l^--\delta^--L^-)>0$, then
it follows that the radius of the free boundary $r^*(h)$ satisfies
\begin{equation}\label{goal}
r^*(h)\approx h^{\frac1{\max(l^\pm-\delta^\pm-1,l^\pm-\delta^\pm-L^\pm)}}.
\end{equation}

We assume now that $A\approx 1$, $B(r)\approx r^m$ and $\Lambda(r)\approx r^{-j}$, where
$m>-1$ and $m+j<1$. After dealing with this situation, we will explain how to deal with the case that $0\le B(r)=O(r^{-1})$ and $\Lambda\approx r^j$ with
$j<2$.
Thus, we assume that there exist positive constants $C_1,C_2$ such that
$C_1\le A\le C_2, C_1r^m\le B(r)\le C_2r^m$ and $C_1r^{-j}\le \Lambda(r)\le C_2r^{-j}$.
From \eqref{f}, it follows that $f'\le0$ and $f''\ge0$. Thus, to get an upper solution, it is sufficient to consider  \eqref{diff-ineq} with
$A=C_2,B=C_2r^m$ and $\Lambda=C_1r^{-j}$.
Similarly, to get a lower solution, it is sufficient to consider \eqref{diff-ineq} with the inequality reversed and with
the above substitutions for $A,B$ and $\Lambda$ with the roles of $C_1$ and $C_2$ reversed.
Substituting in \eqref{diff-ineq} with $f$ as in \eqref{f} and with  $A,B$ and $\Lambda$ as above,
 we will obtain an upper solution if
  the following inequality holds for all large $c$:
\begin{equation}\label{prekeyupper}
\begin{aligned}
&C_2\theta\Big(k(k+1)r^{-k-2}c^{-\delta}+\frac{2lkr^{-k-1}c^{-\delta}}{c-r}+kr^{m-k-1}c^{-\delta}+\frac{l(l-1)c^{-\delta}(\gamma+r^{-k})}{(c-r)^2}+\\
&\frac{r^ml(\gamma+r^{-k})c^{-\delta}}{c-r}\Big)\le \theta^p\thinspace C_1r^{-j}(c-r)^{-l(1-p)}c^{-\delta p}(\gamma+r^{-k})^p,\ R\le r\le c.
\end{aligned}
\end{equation}
Similarly,
we will get a lower solution if the above inequality with the direction of the inequality reversed and with the roles of $C_1$ and $C_2$ reversed
holds for all large $c$.

Multiplying both sides of \eqref{prekeyupper}  by $\theta^{-p}r^{k+2}(c-r)c^\delta$, we rewrite it as
\begin{equation}\label{keyupper}
\begin{aligned}
&C_2\theta^{1-p}\Big(k(k+1)(c-r)+2lkr+kr^{m+1}(c-r)+\frac{l(l-1)\gamma r^{k+2}}{c-r}+\frac{l(l-1)r^2}{c-r}+\\
&l\gamma r^{m+k+2}+lr^{m+2}\Big)
\le C_1r^{k+2-j}(c-r)^{1-l(1-p)}c^{(1-p)\delta}(\gamma+r^{-k})^p,\ R\le r\le c.
\end{aligned}
\end{equation}
Recall that $k=c^{-L}$. If we find values $l=l^+,\delta=\delta^+,L=L^+$, $\gamma=\gamma^+$ and $\theta=\theta^+$, independent of $c$, such that \eqref{keyupper} holds for all large $c$, then we will have an upper solution, while if we find values $l=l^-,\delta=\delta^-, L=L^-$, $\gamma=\gamma^-$ and $\theta=\theta^-$, independent of $c$,
such that \eqref{keyupper} with the roles of $C_1$ and $C_2$ reversed and with the inequality reversed holds for all large $c$, then we will have a lower solution.
Consider for a moment \eqref{keyupper} without the factor $\theta^{1-p}$.
In fact, if we can find values $l^+,\delta^+,L^+,\gamma^+$  such that with $\theta=1$,
for all values of $r\in[R,c]$, the order of the  left hand side of \eqref{keyupper} as $c\to\infty$
is less than or equal  to the order of the right hand side, then we can attain \eqref{keyupper} by making $\theta$ sufficiently
small. Similarly, if we can find values $l^-,\delta^-,L^-,\gamma^-$  such that with $\theta=1$,
for all values of $r\in[R,c]$, the order of the  left hand side of \eqref{keyupper} as $c\to\infty$
is greater than or equal  to the order of the right hand side, then we can attain \eqref{keyupper} with the inequality
reversed and with the roles of $C_1$ and $C_2$ reversed by making $\theta$ sufficiently
large.
Since all the term on the left hand side of \eqref{keyupper} are positive, it suffices to compare the orders of the terms on the
 left hand and right hand
sides, for $r$ in the various appropriate ranges.

Consider  the left hand and right hand sides of \eqref{keyupper} first  for $r=O(1)$ as $c\to\infty$.  The first and third terms on the left hand side
are on the order $c^{1-L}$, while the rest of the terms are of smaller order.
The term on the right hand side is on the order $c^{1-(1-p)l+(1-p)\delta}$.
Thus for an upper solution we need $1-L\le 1-(1-p)l+(1-p)\delta$; that is,
\begin{equation}\label{r=1upper}
l-\delta\le \frac L{1-p}.
\end{equation}
For a lower solution we need
\begin{equation}\label{r=1lower}
l-\delta\ge \frac L{1-p}.
\end{equation}

Now consider $r=\alpha c$ for some $\alpha\in(0,1)$.
The left hand side is on the order $c^{m+2+k}$ and the right hand side is on the order
$c^{(1-p)\delta+1-(1-p)l+k+2-j}$.
Thus, for an upper solution we need
\begin{equation}\label{alpha-c-upper}
l-\delta\le \frac{1-m-j}{1-p},
\end{equation}
while for a lower solution we need
\begin{equation}\label{alpha-c-lower}
l-\delta\ge \frac{1-m-j}{1-p}.
\end{equation}

We have skipped over the transition range between $r=O(1)$ and $r=\alpha c$.
If one substitutes $r=c^\beta$ for some $\beta\in(0,1)$, and makes an analysis of the orders
of the left hand and right hand sides (similar to the analysis we make below in the case
$r=c-c^{-\beta}$), one finds that if \eqref{r=1upper} and \eqref{alpha-c-upper} hold,
then  the order of the left hand side is less than or equal to that of the right hand side, while
if \eqref{r=1lower} and \eqref{alpha-c-lower} hold, then the order on the left hand side is greater
than or equal to that of the right hand side. We leave this to the reader.
Thus, no additional restriction on the parameters is necessary to deal with this range of $r$.

Consider now  the case that $r\to c $. If $1-l(1-p)<-1$, then  the left hand side  of \eqref{keyupper}  blows up at a smaller
order than the right hand side as $r\to c$,
while if  $1-l(1-p)>-1$, then the left hand side blows up at a larger order.
When $1-l(1-p)=-1$, then multiplying both sides by $c-r$ and letting $r\to c$, we find that the resulting limit
on the left hand side will be of smaller or equal order to that of the right hand side if
$(1-p)\delta-j\ge0$.  On the other hand, the resulting limit on the left hand side will be of larger or equal order to that of the
right hand side if  $(1-p)\delta-j\le0$.
Thus, for an upper solution, we need
\begin{equation}\label{r=c-upper}
\begin{aligned}
&l>\frac2{1-p}\ \ \text{or}\\
& l=\frac2{1-p}\ \ \text{and}\ \ \delta\ge \frac j{1-p},
\end{aligned}
\end{equation}
while for a lower solution we need
\begin{equation}\label{r=c-lower}
\begin{aligned}
&l<\frac2{1-p}\ \ \text{or}\\
& l=\frac2{1-p}\ \ \text{and}\ \ \delta\le \frac j{1-p}.
\end{aligned}
\end{equation}

Finally we consider the case that $r$ is close to  $c$. It is enough to consider
$r=c-c^{-\beta}$, with $\beta\ge 0$. The right hand side of \eqref{keyupper} is on the order $c^{k+2-j+\beta\big(l(1-p)-1\big)+(1-p)\delta}$,
 while the left hand side is on the order $c^{\max(k+2+m,k+2+\beta)}$.
Unlike
what occurred in the three ranges of $r$ we  treated above,
here we will need to consider the upper and the lower solution
cases separately.

We first consider the upper solution case.
 When $m>0$ and $\beta\in[0, m)$, then the order of the left hand side of \eqref{keyupper}
 will be less than or equal
to that of the right hand side if $m\le -j+\beta\big(l(1-p)-1\big)+(1-p)\delta$.
For an upper solution,  \eqref{r=c-upper} has already forced us to choose $l\ge\frac2{1-p}$, so we assume here that this condition holds.
Since we require that the above inequality   hold for all such $\beta\in[0,m)$, and
since,  by our assumption, $l(1-p)-1>0$, we obtain the condition
$m\le -j+(1-p)\delta$, or equivalently $\delta\ge\frac{m+j}{1-p}$. (If $m\le0$, then we don't obtain any restriction here.)
Now consider $\beta\ge m\vee\thinspace 0$. Then the order of the left hand side of \eqref{keyupper} will be less than or equal to that of the right hand side if
$\beta\le-j+\beta\big(l(1-p)-1\big)+(1-p)\delta$, or equivalently, if
$\beta\big(2-l(1-p)\big)\le (1-p)\delta-j$.
Since we require this for arbitrarily large $\beta$, we need $l\le\frac2{1-p}$. However, we have already assumed
that $l\ge\frac2{1-p}$. Thus, we conclude that we need $l=\frac2{1-p}$. With this choice, the above inequality becomes
$\delta\ge \frac j{1-p}$. Thus, the analysis for this range of $r$ shows that
in order to have an upper solution, we need
\begin{equation}\label{r^-beta-upper}
l=\frac2{1-p}\ \ \text{and}\ \ \delta\ge \frac{j+m\vee\thinspace 0}{1-p}.
\end{equation}

Now we consider the lower solution case.
When $m>0$ and $\beta\in[0,m)$, then the order of the left hand side of \eqref{keyupper}
will be greater than or equal to that of
the right hand side if  $m\ge -j+\beta\big(l(1-p)-1\big)+(1-p)\delta$.
The above inequality must hold for all $\beta\in[0,m)$. We will assume now that $l=\frac1{1-p}$. (This choice has been made with hindsight.) Note
that this does not interfere with the condition \eqref{r=c-lower} that we already obtained in order to get a lower solution.
With this assumption, the above
inequality will hold for all relevant $\beta$ if $m\ge-j+(1-p)\delta$, or equivalently, if
$\delta\le\frac{m+j}{1-p}$. (If $m\le0$, then we don't obtain any restriction here.)
Now consider $\beta\ge m\vee\thinspace 0$.
 Then the order of the left hand side of \eqref{keyupper} will be greater than or equal to that of
the right hand side if  $\beta\ge -j+\beta\big(l(1-p)-1\big)+(1-p)\delta$. By our choice of $l$, this reduces to
$\beta\ge -j+(1-p)\delta$.  In order for this to hold for all $\beta\ge m\vee\thinspace 0$, we need
$m\vee\thinspace0\ge-j+(1-p)\delta$, or equivalently,
$\delta\le \frac{j+m\vee\thinspace0}{1-p}$. Thus, the analysis of this range of $r$ shows that in order to have a lower solution we can choose
\begin{equation}\label{r^-beta-lower}
l=\frac1{1-p}\ \ \text{and}\ \ \delta\le \frac{j+m\vee\thinspace 0}{1-p}.
\end{equation}

Putting everything together, we see that in order to have an upper solution, it suffices for the parameters to satisfy
\eqref{r=1upper}, \eqref{alpha-c-upper}, \eqref{r=c-upper} and \eqref{r^-beta-upper}.
Recall that we have assumed that $m>-1$ and that  $m+j<1$.
We choose
\begin{equation}\label{upperparameters}
l^+=\frac2{1-p},\ \delta^+=\frac{1+m+j}{1-p},\ L^+=1-m-j.
\end{equation}
In order to have a lower solution, it it suffices for the parameters to satisfy
\eqref{r=1lower}, \eqref{alpha-c-lower}, \eqref{r=c-lower} and \eqref{r^-beta-lower}.
We choose
\begin{equation}\label{lowerparameters}
l^-=\frac1{1-p},\ \delta^-=\frac{m+j}{1-p},\ L^-=1-m-j.
\end{equation}
Note that $l^+-\delta^+=l^--\delta^-$ and $L^+=L^-$. Thus,
recalling the discussion  ending at \eqref{goal},
it follows that
\eqref{goal} holds.
If $m+j\in[0,1)$, then
\begin{equation}\label{final1}
\begin{aligned}
&\max(l^\pm-\delta^\pm-1,l^\pm-\delta^\pm-L^\pm)=l^\pm-\delta^\pm-L^\pm=\frac{1-m-j}{1-p}\thinspace-\thinspace(1-m-j)=\\
&\frac{p(1-m-j)}{1-p},
\end{aligned}
\end{equation}
while if $m+j<0$, then
\begin{equation}\label{final2}
\begin{aligned}
&\max(l^\pm-\delta^\pm-1,l^\pm-\delta^\pm-L^\pm)=l^\pm-\delta^\pm-1=\frac{1-m-j}{1-p}\thinspace-1=\\
&\frac{p-m-j}{1-p}.
\end{aligned}
\end{equation}
Parts (i) and (ii) of the theorem now follow from \eqref{goal}, \eqref{final1} and \eqref{final2}.

We now turn to  the case that
$0\le B(r)= O(r^{-1})$ and $\Lambda\approx r^j$ with
$j<2$. For the upper solution, the worst case for the drift  is $B(r)=\frac{C_2}r$,
 so we will need to consider \eqref{keyupper} with $m=-1$, while for the lower solution, the worst case
 for the drift is $B(r)=0$, so we will need to consider \eqref{keyupper}
(with the inequality reversed)  with the three terms containing $m$ on the left hand side deleted.
For the upper bound everything goes through as above, leading to
the conclusions \eqref{r=1upper}, \eqref{alpha-c-upper}, \eqref{r=c-upper} and \eqref{r^-beta-upper}
 with $m=-1$. For the lower bound,
we obtain
\eqref{r=1lower},  \eqref{r=c-lower} and \eqref{r^-beta-lower} as before.
However, for the case $r=\alpha c$ with $\alpha\in(0,1)$, the left hand side of \eqref{keyupper} is on the order
$c^{1+k}$ instead of on the order $c^{2+m+k}$---that is, it is as if $m=-1$. Consequently, we are lead
to \eqref{alpha-c-lower} with $m$ set at $-1$. Thus, our conclusion is as before, but with $m$ set at $-1$.
This concludes the proof for $r_N^*$.

We now explain how the above calculations also give us the corresponding result for $r_D^*$.
In order to obtain upper and lower solutions for the solution $v_n$ to \eqref{nonLHE-rad-n-Dir},
we use the very same form
for  the test functions $V,f$ as in \eqref{V} and \eqref{f}, and use the very same choices of
$l^\pm,\delta^\pm, L^\pm$. The only difference will be our choice of $c=c(h)$.
Thus, the differential inequalities \eqref{diff-ineq} and \eqref{prekeyupper} will hold when we use $l^+,\delta^+,L^+$, while the reverse
inequalities will hold when we use $l^-,\delta^-,L^-$. Instead of the boundary inequality \eqref{bdryinequ} and its
corresponding reverse inequality,
we need the inequality
\begin{equation}\label{bdryinequ-Dir}
(c-R)^lf(R)\ge h,
\end{equation}
and its corresponding reverse inequality.
Whereas \eqref{bdryinequ} led to \eqref{bdryinequagain},
the inequality \eqref{bdryinequ-Dir} leads to
the inequality
\begin{equation}\label{bdryinequagain-Dir}
\theta(c-R)^lc^{-\delta}(\gamma+R^{-k})\ge h.
\end{equation}
 Whereas  \eqref{bdryinequagain} was fulfilled with
 $c=c(h)\ge Mh^{\frac1{\max(l-\delta-1,l-\delta-L)}}$, for some appropriate $M>0$,
 \eqref{bdryinequagain-Dir} is fulfilled if
 $c=c(h)\ge Mh^\frac1{l-\delta}$, for some appropriate $M>0$, and of course
 the reverse inequality to \eqref{bdryinequagain-Dir} is fulfilled if
 $c=c(h)\le Mh^\frac1{l-\delta}$, for some appropriate $M>0$.
From  \eqref{upperparameters} and \eqref{lowerparameters}, we have
$\frac1{l^+-\delta^+}=\frac1{l^--\delta^-}=\frac{1-p}{1-m-j}$.
Thus,  we obtain $r_D^*(h)\approx h^\frac{1-p}{1-m-j}$.

\hfill $\square$

\medskip

\noindent\it Proof of Theorem \ref{nofreebdry}.\rm\
We first consider   $r_N^*$. Then we show how the proof for $r_N^*$ also works for $r_D^*$.
Virtually all of the work for the proof of this theorem has been done
in the proof of Theorem \ref{outwards}.
In that theorem we thought of $c$ as $c=c(h)$ with $h$ large. In the present case, we don't think of $c$ as a function of $h$.
For the case $m+j>1$,
we want to show that for any given  $h>0$, we can choose $c$ arbitrarily large and find a lower solution $V$ of the form
given in \eqref{V}.   This will show that $r_N^*(h)\ge c$, for any $c$, and thus that $r_N^*(h)=\infty$.
For the case $m+j=1$, we want to show that for sufficiently large $h$, we can choose $c$ arbitrarily large and find a lower solution
$V$ of the form given in \eqref{V}, and we want to show that for sufficiently small $h$, we can choose  $c$  arbitrarily large and find an upper solution $V$ of the form
given in \eqref{V}.
(We need $c$ arbitrarily large even for the upper solution because we only know that the differential inequality
\eqref{keyupper} holds for sufficiently large $c$.)
 This will show that $r_N^*(h)=\infty$ for large $h$ and that $r_N^*(h)<\infty$ for small $h$.
We will prove  part (ii) of the theorem. The small change in the argument that is needed to prove part (i) is
the same as that noted in the final paragraph of the proof of Theorem \ref{outwards} to handle the case that $m=-1$ there.

We first  construct a lower solution in the case that $m+j>1$ or that $m+j=1$ and $h$ is sufficiently large.
The function $V$ in \eqref{V} is constructed from a function $f$ as in \eqref{f}.
The analysis in  the construction of a lower solution in the proof of  Theorem \ref{outwards} goes through verbatim.
We are led to choosing $l^-,\delta^-$ and $L^-$ as in
\eqref{lowerparameters}. This gives $l^-=\frac1{1-p},\ \delta^-=\frac{m+j}{1-p},\ L^-=1-m-j$.
With this choice, for sufficiently large $c$, $V$ will solve the differential inequality \eqref{keyupper} with the direction reversed,
 and with the roles of $C_1$ and $C_2$ reversed, as needed for a lower solution.
We also need $V$ to solve a boundary differential inequality  at $R$.
So when $m+j>1$,  we need
for \eqref{bdryinequagain} to hold with the inequality reversed, for any fixed $h$, and for arbitrarily large $c$, while
for $m+j=1$, we need for this inequality to hold for sufficiently large $h$, and for arbitrarily large $c$.
Substituting in \eqref{bdryinequagain} the values obtained above for $l,\delta$ and $k=c^{-L}$, and reversing the inequality, we obtain the inequality
\begin{equation}\label{bdryinequyetagain}
\theta\big[\frac1{1-p}(c-R)^{\frac p{1-p}}\thinspace c^{-\frac{m+j}{1-p}}(\gamma+R^{-c^{m+j-1}})+
(c-R)^\frac1{1-p}\thinspace c^{-\frac{m+j}{1-p}}c^{m+j-1}R^{-c^{m+j-1}-1}\big]\le h.
\end{equation}
When $m+j>1$, the left hand side of
 \eqref{bdryinequyetagain} converges to 0 when $c\to\infty$; thus, indeed  \eqref{bdryinequyetagain} holds for any $h>0$ and arbitrarily large
$c$. When $m+j=1$, the left hand side of  \eqref{bdryinequyetagain} converges to $\theta R^{-2}$ when $c\to\infty$; thus,  \eqref{bdryinequyetagain} holds
for sufficiently large $h$ and arbitrarily large $c$.

We now turn to the construction of an upper solution in the case that $m+j=1$ and $h$ is sufficiently small.
The analysis in  the construction of an upper solution in the proof of  Theorem \ref{outwards} goes through verbatim.
We are led to choosing $l^+,\delta^+$ and $L^+$ as in
\eqref{upperparameters}. This gives $l^+=\frac2{1-p},\ \delta^+=\frac{1+m+j}{1-p},\ L^+=1-m-j$.
With this choice, for sufficiently large $c$,  $V$ will solve the differential inequality \eqref{keyupper} as needed for an upper solution.
We also need $V$ to solve a boundary differential inequality  at $R$.
So   we need
for \eqref{bdryinequagain} to hold for sufficiently small $h$ and arbitrarily large  $c$.
Substituting in \eqref{bdryinequagain} the values obtained above for $l,\delta$ and $k=c^{-L}$, and using the fact that $m+j=1$, we obtain the inequality
\begin{equation}\label{bdryinequyetagainagain}
\theta\big[\frac1{1-p}(c-R)^{\frac p{1-p}}\thinspace c^{-\frac1{1-p}}(\gamma+R^{-1})+
(c-R)^\frac1{1-p}\thinspace c^{-\frac1{1-p}}R^{-2}\big]\ge h.
\end{equation}
The left hand side of \eqref{bdryinequyetagainagain} converges to $\theta R^{-2}$ when $c\to\infty$; thus
\eqref{bdryinequyetagainagain} holds for sufficiently small $h$ and arbitrarily large $c$.
This completes the proof for $r_N^*$.

For $r_D^*$, the only change is that the boundary inequality required is different. For the case of $m+j>1$ or $m+j=1$ and $h$ sufficiently
large, instead of the boundary inequality
\eqref{bdryinequagain} with the inequality reversed, we need
\begin{equation*}
\theta(c-R)^lc^{-\delta}(\gamma+R^{-k})\le h.
\end{equation*}
Substituting the values appearing  in the penultimate  paragraph above  for $l^-,\delta^-$ and $k=c^{-L^-}$, we obtain the inequality
\begin{equation}\label{bdryinequyetagainDir}
\theta(c-R)^\frac1{1-p}c^{-\frac{m+j}{1-p}}(\gamma+R^{-c^{m+j-1}})\le h.
\end{equation}
When $m+j>1$, the left hand side of
 \eqref{bdryinequyetagainDir} converges to 0 when $c\to\infty$; thus, indeed  \eqref{bdryinequyetagainDir} holds for any $h>0$ and arbitrarily large
$c$. When $m+j=1$, the left hand side of  \eqref{bdryinequyetagainDir} converges to $\theta(\gamma+ R^{-1})$; thus,  \eqref{bdryinequyetagainDir} holds
for sufficiently large $h$ and arbitrarily large $c$.
For the case $m+j=1$ and $h$ sufficiently small, we need
\eqref{bdryinequyetagainDir} with the inequality reversed; that is,
\begin{equation}\label{bdryinequyetagainagainDir}
\theta(c-R)^\frac1{1-p}c^{-\frac1{1-p}}(\gamma+R^{-1})\ge h.
\end{equation}
The left hand side of \eqref{bdryinequyetagainagainDir} converges to $\theta(\gamma+R^{-1})$ when $c\to\infty$; thus,
\eqref{bdryinequyetagainagainDir}  holds for sufficiently small $h$ and arbitrarily large $c$.
This completes the proof for $r_D^*$.
\hfill $\square$

\section{Proof of Theorem \ref{m>-1}}\label{logh}
As with the proofs from the last section, we present the proof for $r_N^*$, and then,  upon its completion,  we explain how the proof for $r_D^*$ follows from the
proof for $r_N^*$. For the solution $u_n$ to \eqref{nonLHE-rad-n},
we look for upper and lower solutions
  $V$ in the form \eqref{V} with $l=\frac2{1-p}$.
In order for $V$ to be an upper  solution,
the function $f$ appearing in the definition of $V$ must satisfy \eqref{diff-ineq} and \eqref{bdryinequ} with
$l=\frac2{1-p}$.
In order for $V$ to be a lower solution, the function $f$ must satisfy  \eqref{diff-ineq} and \eqref{bdryinequ} with
$l=\frac2{1-p}$ and with the inequality reversed.
If such a $V$ with $c=c(h)$ is an upper (lower) solution, then it follows that $r^*(h)\le c(h)$ $\big(r^*(h)\ge c(h)\big)$.

For the lower solution, we consider $f$ in the form
\begin{equation}\label{fagain}
f=e^{2kc^{m+1}-kr^{m+1}},
\end{equation}
with $k>0$ independent of $c$.
Then \eqref{bdryinequ} with  $l=\frac2{1-p}$ and with the sign of the inequality reversed becomes
\begin{equation}\label{bdryinequ-m>-1low}
(c-R)^\frac2{1-p}k(m+1)R^me^{2 kc^{m+1}-kR^{m+1}}+\frac2{1-p}(c-R)^\frac{1+p}{1-p}e^{2 kc^{m+1}-kR^{m+1}}\le h.
\end{equation}
For sufficiently large $c$, the left hand side of \eqref{bdryinequ-m>-1low} is bounded from above by $e^{(2k+1)c^{m+1}}$.
Thus, for sufficiently large $h$, \eqref{bdryinequ-m>-1low} will hold if we choose
$c=c(h)=C_0(\log h)^{\frac1{m+1}}$, for an appropriate $C_0>0$.
We will now show that we can pick  $k$ in the definition of $f$ so that
\eqref{diff-ineq} holds with $l=\frac2{1-p}$ and  with the sign of the inequality reversed.
It will then follow that $V$ is a lower solution, and consequently   that
$r^*(h)\ge c(h)= C_0(\log h)^{\frac1{m+1}}$, for large $h$, thereby proving the lower bound in the theorem.

We have
$f'=-k(m+1)r^mf$ and $f''=k^2(m+1)^2r^{2m}f-km(m+1)r^{m-1}f$.
By assumption, there exist constants $C_1,C_2>0$ and $N>0$ such that
$C_1\le A(r)\le C_2,-C_2r^m\le B(r)\le -C_1r^m,C_1r^{-N}\le\Lambda(r)\le C_2r^N$.
Thus, in order to get a lower solution, in \eqref{diff-ineq},  it suffices to consider the reverse inequality and to substitute
$-C_2r^m$ for $B$ (since $f'\le0$), to substitute $C_2r^N$ for $\Lambda$, to substitute
$C_1k^2(m+1)^2r^{2m}f-C_2km(m+1)r^{m-1}f$  for $Af''$,  to substitute $C_1f$ for $Af$
and to substitute $C_1f'$ for $Af'$ (since $f'\le0$).
Recall also that $l=\frac2{1-p}$.
Making these substitutions in \eqref{diff-ineq},
 multiplying both sides by $\frac1f$ and grouping certain terms, we obtain
\begin{equation}\label{keym>-1}
\begin{aligned}
&k(m+1)\Big(k(m+1)C_1-C_2\Big)r^{2m}+\frac{4C_1k(m+1)-2C_2}{(1-p)(c-r)}r^m-
C_2km(m+1)r^{m-1}+\\
&\frac{2C_1(1-p)}{(1-p)^2(c-r)^2}\ge\frac{C_2r^N}{(c-r)^2}
e^{-(1-p)\big( 2kc^{m+1}-kr^{m+1}\big)}, r\le R\le c.
\end{aligned}
\end{equation}
The maximum of $C_2 r^Ne^{-(1-p)\big(2 kc^{m+1}-kr^{m+1}\big)}$ over $r\in[R,c]$,
occurs of course at $r=c$ and is equal to $C_2c^Ne^{-(1-p)k c^{m+1}}$. Thus, for any $k>0$, it follows that for all large $c$,
this last term is no greater than $\frac{2C_1(1+p)}{(1-p)^2}$.
Thus, in order to obtain \eqref{keym>-1}, it suffices to choose $k>0$ so  that
$4C_1k(m+1)-2C_2\ge0$ and
$$
k(m+1)\Big(k(m+1)C_1-C_2\Big)r^{2m}\ge
C_2km(m+1)r^{m-1}, r\le R\le c.
$$
Since $m>-1$, clearly the first inequality holds for all sufficiently large $k$.
Since $m>-1$, we have $2m>m-1$ and thus the second inequality also holds for sufficiently large $k$.

We now turn to the upper solution.
We consider $f$ in the form
$$
f=\gamma e^{\frac12kc^{m+1}-kr^{m+1}},
$$
with $k,\gamma>0$, independent of $c$.
Then \eqref{bdryinequ} with  $l=\frac2{1-p}$ becomes
\begin{equation}\label{bdryinequ-m>-1up}
(c-R)^\frac2{1-p}\gamma k(m+1)R^me^{\frac12 kc^{m+1}-kR^{m+1}}+\frac2{1-p}(c-R)^\frac{1+p}{1-p}\gamma e^{\frac12 kc^{m+1}-kR^{m+1}}\ge h.
\end{equation}
For sufficiently large $c$ (depending on $k$ and $\gamma$), the left hand side of \eqref{bdryinequ-m>-1up} is bounded from below by $e^{\frac12kc^{m+1}}$.
Thus, for sufficiently large $h$ (depending on $k$ and $\gamma$), \eqref{bdryinequ-m>-1up} will hold if we choose
$c=c(h)=C_0(\log h)^{\frac1{m+1}}$, for an appropriate  $C_0>0$.
We will now show that we can pick  $k,\gamma>0$ in the definition of $f$ so that
\eqref{diff-ineq} holds with $l=\frac2{1-p}$.
It will then follow that $V$ is an upper solution, and consequently   that
$r^*(h)\ge c(h)= C_0(\log h)^{\frac1{m+1}}$, for large $h$, thereby proving the upper bound in the theorem.

The same type of  considerations that led to \eqref{keym>-1} show that in order to obtain \eqref{diff-ineq}, it suffices to
obtain the inequality
 \begin{equation}\label{keym>-1upper}
\begin{aligned}
&k(m+1)\Big(k(m+1)C_2-C_1\Big)r^{2m}+\frac{4C_2k(m+1)-2C_1}{(1-p)(c-r)}r^m-
C_1km(m+1)r^{m-1}+\\
&\frac{2C_2(1-p)}{(1-p)^2(c-r)^2}\le\gamma^{p-1}\frac{C_1r^{-N}}{(c-r)^2}
e^{-(1-p)\big( \frac12kc^{m+1}-kr^{m+1}\big)}, R\le r\le c.
\end{aligned}
\end{equation}
Now choose $k>0$ sufficiently small  so that $k(m+1)C_2-C_1<0$ and $4C_2k(m+1)-2C_1<0$.
Then the first three terms on the left hand side of \eqref{keym>-1upper} are negative.
For say, $r\ge(\frac34)^\frac1{m+1}c$, and $c$ sufficiently large,  the last term on the left hand side of \eqref{keym>-1upper} is smaller
than the right hand side of \eqref{keym>-1upper}.
Thus, we conclude that the inequality in \eqref{keym>-1upper} holds for $(\frac34)^\frac1{m+1}c\le r\le c$.
On the other hand, for $R\le r\le (\frac34)^\frac1{m+1}c$, the first and third terms on the left hand side
of \eqref{keym>-1upper} are negative and bounded from 0, the second term is negative and the last term
is on the order of $\frac1{c^2}$. Thus, for sufficiently large $c$, the left hand side
of \eqref{keym>-1upper} is negative for   $R\le r\le (\frac34)^\frac1{m+1}c$.
Thus, we conclude that
\eqref{keym>-1upper} holds for the entire range of $r$ as specified.
This completes the proof for $r_N^*$.

We now explain how the above calculations also give us the corresponding result for $r_D^*$.
In order to obtain upper and lower solutions for the solution $v_n$ to \eqref{nonLHE-rad-n-Dir},
we use the very same form
as used above, the only difference being that
we consider varying  $c=c(h)$.  For the lower solution obtained above for \eqref{nonLHE-rad-n}, the boundary inequality took on the form \eqref{bdryinequ-m>-1low}.
For the lower solution for \eqref{nonLHE-rad-n-Dir}, the boundary inequality takes on the form
\begin{equation}\label{Dirbclow}
(c-R)^\frac2{1-p}\thinspace e^{2kc^{m+1}-kR^{m+1}}\le h.
\end{equation}
For the upper solution obtained above for \eqref{nonLHE-rad-n}, the boundary inequality took on the form \eqref{bdryinequ-m>-1up}.
For the upper solution for \eqref{nonLHE-rad-n-Dir}, the boundary inequality takes on the form
\begin{equation}\label{Dirbcup}
(c-R)^\frac2{1-p}\thinspace\gamma e^{\frac12kc^{m+1}-kR^{m+1}}\ge h.
\end{equation}
It is easy to see that just as with \eqref{bdryinequ-m>-1low} and \eqref{bdryinequ-m>-1up}, we
can fulfill \eqref{Dirbclow} and \eqref{Dirbcup} by choosing
$c=c(h)=C_0(\log h)^{\frac1{m+1}}$, for an appropriate constant $C_0>0$ in each case. Thus,
 as with $r_N^*$, we obtain, $r_D^*(h)\approx (\log h)^\frac1{m+1}$.
\hfill $\square$

\section{Proof of  Theorem \ref{-1}}\label{fuchsian}
As with the proofs from the previous two sections, we present the proof for $r_N^*$, and then, upon its completion, we
explain how the proof for $r_D^*$ follows from the proof for $r_N^*$.
For the solution $u_n$ to \eqref{nonLHE-rad-n}, we look for  upper and lower solutions  $V$ in the form \eqref{V} with $l=\frac2{1-p}$. In order for $V$ to be an upper  solution,
the function $f$ appearing in the definition of $V$ must satisfy \eqref{diff-ineq} and \eqref{bdryinequ} with
$l=\frac2{1-p}$.
In order for $V$ to be a lower solution, the function $f$ must satisfy  \eqref{diff-ineq} and \eqref{bdryinequ} with
$l=\frac2{1-p}$ and with the inequality reversed.
If such a $V$ with $c=c(h)$ is an upper (lower) solution, then it follows that $r^*(h)\le c(h)$ $\big(r^*(h)\ge c(h)\big)$.

Recalling that $A=A_0$, that $B=-\frac{B_0}r$ and that $\mu=\frac{B_0}{A_0}>0$, we can write \eqref{diff-ineq} with $l=\frac2{1-p}$
as
\begin{equation}\label{originequ-f}
f''+\big(\frac\mu r-\frac4{(1-p)(c-r)}\big)f'+\big(\frac{2(1+p)}{(1-p)^2(c-r)^2}-\frac{2\mu}{(1-p)r(c-r)}\big)f\le \frac{\Lambda f^p}{A_0(c-r)^2}.
\end{equation}

Consider the solution
$g=g_c>0$ to the linear equation
\begin{equation}\label{originequ-g}
\begin{aligned}
&g''+\big(\frac\mu r-\frac4{(1-p)(c-r)}\big)g'-\frac{2\mu}{(1-p)r(c-r)}g= 0, \ R\le r<c;\\
&g(R)=1;\\
&g\le1\ \text{is maximal}.
\end{aligned}
\end{equation}
By maximal, we mean that $g_c=\lim_{n\to\infty}g_{c;n}$, where $g_{c;n}$ satisfies
\begin{equation}\label{originequ-gn}
\begin{aligned}
&g''+\big(\frac\mu r-\frac4{(1-p)(c-r)}\big)g'-\frac{2\mu}{(1-p)r(c-r)}g= 0, \ R\le r\le c-\frac1n;\\
&g(R)=g(c-\frac1n)=1.
\end{aligned}
\end{equation}
(Note that by the maximum principle, $g_{c;n}\le1$.)
We will show that
\begin{equation}\label{gboundedfrom0}
g_c\ \text{is  bounded from}\  0,\  \text{uniformly in}\ c.
\end{equation}
Thus, recalling that $\Lambda$ is assumed to be bounded and bounded from 0,
 it follows that by  choosing $\delta>0$ sufficiently small (independent of $c$), the function $f^+_c=\delta g_c$ will satisfy \eqref{originequ-f}
and the function $f_c^-=\delta^{-1} g_c$ will satisfy \eqref{originequ-f} with the inequality reversed.
The boundary flux condition \eqref{bdryinequ} for an upper solution then becomes
\begin{equation}\label{g'R+}
\frac2{1-p}(c-R)^{\frac{1+p}{1-p}}\delta-(c-R)^\frac2{1-p}\delta g'(R)\ge h,
\end{equation}
and for a lower solution it becomes
\begin{equation}\label{g'R-}
\frac2{1-p}(c-R)^{\frac{1+p}{1-p}}\delta^{-1}-(c-R)^\frac2{1-p}\delta^{-1} g'(R)\le h.
\end{equation}
We will prove that
\begin{equation}\label{g'bound}
C_1c^{\mu-1}\le -g_c'(R)\le C_2c^{\mu-1}, \ \text{where}\ C_1,C_2>0 \ \text{are independent of}\ c.
\end{equation}
From \eqref{g'R+}-\eqref{g'bound} and the fact that $\mu>0$, it follows that we can obtain an upper solution for large $h$
by choosing an appropriate $c=c(h)$ with $c^{\frac2{1-p}+\mu-1}$ on the order $h$,
and also a lower solution for large $h$ by choosing
an appropriate $c=c(h)$  with $c^{\frac2{1-p}+\mu-1}$  on the order $h$; that is, for both the upper and the lower
solution we  have $c(h)$ on the order $h^{\frac{1-p}{1+p+\mu(1-p)}}$.
Thus, we conclude that $r^*(h)\approx h^{\frac{1-p}{1+p+\mu(1-p)}}$.
Therefore, to conclude the proof, we need to prove \eqref{gboundedfrom0} and \eqref{g'bound}.

To prove  \eqref{gboundedfrom0}, we construct a lower solution for \eqref{originequ-g} that is bounded from 0, uniformly over $c$.
Let $v_c$ satisfy
\begin{equation}\label{expectedhitting}
\begin{aligned}
&v''+\big(\frac\mu r-\frac4{(1-p)(c-r)}\big)v'=-\frac1{r(c-r)},\ R\le r<c;\\
& v(R)=0;\\
&v\ge 0\ \text{is minimal}.
\end{aligned}
\end{equation}
By minimal, we mean that $v_c=\lim_{n\to\infty}v_{c;n}$, where
$v_{c;n}$ solves
\begin{equation}\label{expectedhittingn}
\begin{aligned}
&v''+\big(\frac\mu r-\frac4{(1-p)(c-r)}\big)v'=-\frac1{r(c-r)},\ R\le r\le c-\frac1n;\\
& v(R)=v(c-\frac1n)=0.
\end{aligned}
\end{equation}
Solving  \eqref{expectedhittingn} gives
$$
v_{n;c}(r)=k_{n;c}\int_R^rs^{-\mu}(c-s)^{-\frac4{1-p}}ds-\int_R^rdz\thinspace z^{-\mu}(c-z)^{-\frac4{1-p}}\int_R^zs^{\mu-1}(c-s)^{\frac{3+p}{1-p}}ds,
$$
where
$$
k_{n;c}=\frac{\int_R^{c-\frac1n}dz\thinspace     z^{-\mu}(c-z)^{-\frac4{1-p}}\int_R^zs^{\mu-1}(c-s)^{\frac{3+p}{1-p}}ds}{\int_R^{c-\frac1n}s^{-\mu}(c-s)^{-\frac4{1-p}}ds}.
$$
Since $\lim_{n\to\infty}\int_R^{c-\frac1n}s^{-\mu}(c-s)^{-\frac4{1-p}}ds=\infty$, we have
$$\lim_{n\to\infty}k_{n;c}=\int_R^cs^{\mu-1}(c-s)^{\frac{3+p}{1-p}}ds.
$$
Thus we obtain
\begin{equation}\label{expectedhitting-v}
v_c(r)=\int_R^rdz\thinspace z^{-\mu}(c-z)^{-\frac4{1-p}} \int_z^cs^{\mu-1}(c-s)^{\frac{3+p}{1-p}}ds.
\end{equation}
It is easy to check that $e^{-\frac{2\mu}{1-p}v_{n;c}}$ is a lower solution for \eqref{originequ-gn}; it satisfies
the boundary condition in \eqref{originequ-gn} and it satisfies the differential inequality obtained by changing
the equal sign to $\ge$.
Thus, by the maximum principle, $g_{n;c}\ge e^{-\frac{2\mu}{1-p}v_{n;c}}$, and consequently
\begin{equation}\label{gcvc}
g_c\ge e^{-\frac{2\mu}{1-p}v_c}.
\end{equation}
To complete the proof of  \eqref{gboundedfrom0} we now show that    $\sup_{c>R}v_c(c)<\infty$.

Making
the change of variables $s=tc$ in the inside integral, and then making the change of variables $z=yc$ in the outside
integral, we obtain
\begin{equation}
v_c(c)=\int_{\frac Rc}^1dy\thinspace (1-y)^{-\frac4{1-p}}\int_y^1t^{\mu-1}(1-t)^{\frac{3+p}{1-p}}dt
\end{equation}
Thus, $\lim_{c\to\infty}v_c(c)=\int_0^1dy\thinspace (1-y)^{-\frac4{1-p}}\int_y^1t^{\mu-1}(1-t)^{\frac{3+p}{1-p}}dt<\infty$.
This completes the proof of \eqref{gboundedfrom0}.

We now prove \eqref{g'bound}. Writing the differential equation in \eqref{originequ-g} as
$g''+\big(\frac\mu r-\frac4{(1-p)(c-r)}\big)g'=\frac{2\mu}{(1-p)r(c-r)}g$, and treating the right hand side as an inhomogeneous term,
we solve  by the method used above for $v_c$, obtaining the integral equation
\begin{equation}\label{inteq}
g_c(r)=1-\int_R^rdz\thinspace z^{-\mu}(c-z)^{-\frac4{1-p}} \int_z^cs^{\mu-1}(c-s)^{\frac{3+p}{1-p}}g_c(s)ds.
\end{equation}
By \eqref{gboundedfrom0}, we can choose $\gamma>0$, independent of $c$,  such that $\gamma\le \inf_{R\le r<c}g_c(r)$. Also, of course, $g_c\le 1$.  Thus, from \eqref{inteq} we have
\begin{equation}\label{upperlower}
\begin{aligned}
&-g_c'(R)=\lim_{r\to R^+}\frac{1-g_c(r)}{r-R}\ge\gamma R^{-\mu}(c-R)^{-\frac4{1-p}} \int_R^cs^{\mu-1}(c-s)^{\frac{3+p}{1-p}}ds;\\
&-g_c'(R)=\lim_{r\to R^+}\frac{1-g_c(r)}{r-R}\le R^{-\mu}(c-R)^{-\frac4{1-p}} \int_R^cs^{\mu-1}(c-s)^{\frac{3+p}{1-p}}ds.
\end{aligned}
\end{equation}
We have
\begin{equation}\label{finalest}
\begin{aligned}
&\int_R^cs^{\mu-1}(c-s)^{\frac{3+p}{1-p}}ds\le(c-R)^{\frac{3+p}{1-p}}\int_R^cs^{\mu-1}ds=
(c-R)^{\frac{3+p}{1-p}}\thinspace\frac{c^\mu-R^\mu}\mu;\\
&\int_R^cs^{\mu-1}(c-s)^{\frac{3+p}{1-p}}ds\ge\int_{\frac c4}^{\frac c2}s^{\mu-1}(c-s)^{\frac{3+p}{1-p}}ds\ge
\frac14\frac14^{\mu-1}\frac12^{\frac{3+p}{1-p}}c^{\frac{3+p}{1-p}+\mu}.
\end{aligned}
\end{equation}
Now \eqref{g'bound} follows from \eqref{upperlower} and \eqref{finalest}.
This completes the proof for $r_N^*$.

We now explain how the above calculations also give us the corresponding result for $r_D^*$.
In order to obtain upper and lower solutions for the solution $v_n$ to \eqref{nonLHE-rad-n-Dir},
we use the very same form
as used above, the only difference being that
we  vary  $c=c(h)$.
 For the upper solution obtained above for \eqref{nonLHE-rad-n}, the boundary inequality took on the form \eqref{g'R+}.
For the upper solution for \eqref{nonLHE-rad-n-Dir}, the boundary inequality takes on the form
\begin{equation}\label{Dirfuchup}
(c-R)^\frac2{1-p}\delta\ge h.
\end{equation}
 For the lower solution obtained above for \eqref{nonLHE-rad-n}, the boundary inequality took on the form \eqref{g'R-}.
For the lower solution for \eqref{nonLHE-rad-n-Dir}, the boundary inequality takes on the form
\begin{equation}\label{Dirfuchlow}
(c-R)^\frac2{1-p}\delta^{-1}\le h.
\end{equation}
Clearly, \eqref{Dirfuchup} and \eqref{Dirfuchlow} can each be fulfilled by choosing
$c=c(h)=C_0h^\frac{1-p}2$, for an appropriate constant $C_0>0$ in each case. Thus, we obtain
$r_D^*(h)\approx h^\frac{1-p}2$.
\hfill$\square$

\section{Proofs of Theorems \ref{Laplacian} and \ref{nonradial}}\label{generalcase}
We will give the proof of the two theorems for the case of $r_N^*$. As will be clear from the proof below,
the proof for $r_D^*$ is virtually the same, except
that we don't need to assume that $D$ is star-shaped because now the boundary condition does not contain a derivative.

To prove the two theorems, we simply sketch how to
 convert the upper and lower solutions obtained in the radial case to upper and lower solutions
in the non-radial case. Since $0\in D$, there exist  $0<R<R_1$ such that
$B_R\subset D\subset B_{R_1}$.  Let $V(r)$ be a radial function defined on $R^d-B_R$.
For  $x\in\partial D$, we have $(\nabla V\cdot\bar n)(x)=\gamma(x)V'(r)|_{r=|x|}$, for some
function $\gamma$ defined on $\partial D$. Of course, in general, $-1\le \gamma\le 1$; however, by the
assumption that $D$ is star-shaped with respect to 0, it follows that there exists a $\delta>0$ such that
$\delta\le\gamma\le 1$.
Let $\{V=V_c=V_{c(h)}\}$, for large $h$,  be one of the collections of  upper (lower) solutions that we constructed for the proofs of the other theorems.
They were all sufficiently regular so that if $V'(R)\le -h$ ($V'(R)\ge -h$), then there are  positive constants
$\gamma_1<1$ and $\gamma_2>1$ such that $V'(r)\le-\gamma_1 h$ ($V'(r)\ge -\gamma_2 h$), for all $r\in[R,R_1]$.
Thus from the point of view of the boundary flux condition, the collection $\{V_{c(\frac h{\gamma_1})}\}$
($\{V_{c(\frac h{\gamma_2})}\}$) would satisfy the upper solution (lower solution) boundary flux inequality on $\partial D$.
This collection of course  gives the same asymptotic order for $r^*(h)$ as does
the collection $\{V_{c(h)}\}$.

In light of Theorem \ref{-1}, the above paragraph is enough to  complete the proof of Theorem \ref{Laplacian}.
For Theorem \ref{nonradial}, we also
need to  consider the differential inequalities inside the domain ($R^d-\bar B_R$ or $R^d-\bar D$) for upper and lower solutions.

Consider the operator $L$ as in \eqref{nondivL}. Recall that by the assumption in Theorem \ref{nonradial}, this operator
is defined on all of $R^d$; thus in particular, it is defined on $R^d-B_R$.
Let
$$
A(x)=\sum_{i,j=1}^da_{i,j}(x)\frac{x_i}{|x|}\frac{x_j}{|x|}.
$$
For a function $V$ depending only on $r=|x|$, we   have
$$
LV(x)=A(x)V''(r)
-\frac{\sum_{i=1}^da_{i,i}(x)-A(x)}rV'(r)-\frac xr\cdot b(x)V'(r)
$$
By \eqref{unifell}, we have $C_1\le A(x)\le C_2$ and
$0\le \sum_{i=1}^da_{i,i}(x)\le C_2d$.
Note that
\begin{equation}\label{Ainequ}
\begin{aligned}
& C_1V''(r)\le A(x)V''(r)\le C_2V''(r),\ \text{if}\ V''(r)\ge0;\\
&C_2V''(r)\le A(x)V''(r)\le C_1V''(r),\ \text{if}\ V''(r)<0.
\end{aligned}
\end{equation}
Assume now that $V'\le0$, and note that all of the  radially symmetric upper and lower solutions
we constructed satisfy this condition.
Recalling the definition  \eqref{radialb}, it follows that if
$B^+(r)\le K r^m$,  for some $m>-1$ and some $K\in R$,  then for any $\epsilon>0$, for sufficiently large $r$, one has
$$
-\frac{\sum_{i=1}^da_{i,i}(x)-A(x)}rV'(r)-\frac xr\cdot b(x)V'(r)\le (K+\epsilon)r^m|V'(r)|,
$$
while if $B^-(r)\ge Kr^m$, then for any $\epsilon>0$, for sufficiently large $r$, one has
$$
-\frac{\sum_{i=1}^da_{i,i}(x)-A(x)}rV'(r)-\frac xr\cdot b(x)V'(r)\ge (K-\epsilon)r^m|V'(r)|.
$$
Finally, recalling \eqref{AB},  note that if $AB^-\ge0$, then
 $$
 -\frac{\sum_{i=1}^da_{i,i}(x)-A(x)}rV'(r)-\frac xr\cdot b(x)V'(r)\ge0.
 $$
Using the above facts, it is easy to verify that for each of the  various parts of Theorem \ref{nonradial}, the
condition given there is
 enough to guarantee that the  radial upper and lower solutions constructed for the corresponding Theorem mentioned there work as upper
and lower solutions for the non-radial case.
\hfill $\square$


\begin{thebibliography}{99}



\bibitem{BS85}  Bandle, C. and Stakgold, I., \emph{Reaction-diffusion and dead cores}, Free boundary problems: application and theory, Vol. IV, 436-448, Res. Notes in Math., 121, Pitman, Boston, MA, 1985.


\bibitem{P88} Pinsky, R. \emph{The dead core for reaction-diffusion equations with convection and its connection with the first exit time of the related Markov diffusion process}, Nonlinear Anal. 12 (1988), 451-471.

\bibitem{P92}  Pinsky, R. \emph{The interplay of nonlinear reaction and convection in dead core behavior for reaction-diffusion equations},
 Nonlinear Anal. 18 (1992), 1113-1123.

\bibitem{P14} Pinsky, R. \emph{Universal bound independent of geometry for solution to symmetric diffusion equation in exterior domain with boundary flux},
preprint

\bibitem{S73} Sattinger, D., \emph{Topics in stability and bifurcation theory}, Lecture Notes in Mathematics, Vol. 309, Springer-Verlag, Berlin-New York, 1973.

\bibitem{S96} Sperb, R. \emph{Some complementary estimates in the dead core problem}, Nonlinear problems in applied mathematics, 217-224, SIAM, Philadelphia (1996).

\end{thebibliography}
\end{document}